\documentclass[11pt]{amsart}
\usepackage{graphicx}
\usepackage{amssymb}
\usepackage{fullpage}
\usepackage{natbib}
\usepackage{algorithm}
\usepackage{algpseudocode}
\usepackage{placeins}
\usepackage{xcolor}
\usepackage[unicode]{hyperref}
\usepackage{amsmath,amsthm,enumerate,bbm,tabularx,multicol,environ}
\usepackage[nolist]{acronym}
\usepackage{threeparttable,booktabs}

\newcommand{\ba}{\boldsymbol{a}}

\newcommand{\bd}{\boldsymbol{d}}
\newcommand{\be}{\boldsymbol{e}}
\newcommand{\bg}{\boldsymbol{g}}

\newcommand{\bx}{\boldsymbol{x}}
\newcommand{\bu}{\boldsymbol{u}}
\newcommand{\by}{\boldsymbol{y}}

\newcommand{\bxi}{\boldsymbol{\xi}}

\newcommand{\bz}{\boldsymbol{z}}

\newcommand{\bv}{\boldsymbol{v}}
\newcommand{\bzero}{\boldsymbol{0}}

\newcommand{\bnu}{\boldsymbol{\nu}}

\newcommand{\btheta}{\boldsymbol{\theta}}

\DeclareMathOperator*{\argmin}{argmin}
\DeclareMathOperator*{\argmax}{argmax}

\newcommand{\field}[1]{\mathbb{#1}}

\newcommand{\R}{\field{R}}

\newcommand{\E}{\field{E}}

\newcommand{\cV}{\mathcal{V}}
\newcommand{\cCV}{\mathcal{C}_{\mathcal{V}}}

\DeclareMathOperator{\Regret}{Regret}

\newcommand{\indicator}{\mathbf{1}}

\DeclareFontFamily{U}{mathx}{\hyphenchar\font45}
\DeclareFontShape{U}{mathx}{m}{n}{<-> mathx10}{}
\DeclareSymbolFont{liftaccents}{U}{mathx}{m}{n}
\DeclareMathAccent{\widecheck}{\mathord}{liftaccents}{"71}

\newcommand{\decode}{{\sf decode}}

\newtheorem{theorem}{Theorem}[section]
\newtheorem{proposition}[theorem]{Proposition}

\newtheorem{lemma}[theorem]{Lemma}
\theoremstyle{definition}

\newcounter{assump}
\theoremstyle{remark}
\newtheorem{remark}[theorem]{Remark}
\usepackage{tikz}
\usetikzlibrary{arrows.meta}
\usepackage{mathtools}
\mathtoolsset{showonlyrefs}
\usepackage{amsaddr}

\hypersetup{
    colorlinks=true,
    linkcolor=blue,
    citecolor=blue,
    urlcolor=blue
}

\makeatletter
\AtBeginDocument{%
  \let\saved@float@makebox\float@makebox
  \renewcommand{\float@makebox}[1]{%
    \color@vbox\normalcolor
      \saved@float@makebox{#1}%
    \color@endbox
  }%
}
\makeatother
\usepackage{iroshizuku-hyperref}
\SetIroTheme{midoriao}

\author{Julien Grand-Cl\'{e}ment}
\address{ISOM Department, HEC Paris\\
  78350, Jouy-en-Josas, France}
\email{grand-clement@hec.fr}
\author{Joon Kwon}
\address{INRAE \& Université Paris--Saclay\\91120, Palaiseau, France}
\email{joon.kwon@inrae.fr}
\author{Francesco Orabona}
\address{King Abdullah University of Science and Technology (KAUST)\\
23955--6900, Thuwal, Kingdom of Saudi Arabia}\email{francesco@orabona.com}
\date{\today}
\title{Conic-Lifting for Learning-rate-free Convex Optimization}

\begin{document}

\begin{abstract}
A central topic in online and stochastic optimization is to design algorithms that do not require tuning of hyperparameters. Successful strategies in online learning include AdaGrad-style learning rates, Conic Blackwell algorithms, and parameter-free algorithms. However, these classes require different, often difficult, analyses.
We develop a unifying family of \ac{OCO} algorithms that require no user-tuned learning rate on convex domains that may be unbounded. The construction lifts the domain to a cone by adding one scalar coordinate, an idea going back to \citet{blackwell1954controlled}. The lifted subgradients are orthogonal to the corresponding lifted decisions, leaving the comparator scale free to optimize in the regret analysis.
Applying \ac{FTRL} and \ac{OMD} on the cone recovers \ac{CBA} and \ac{CBA+} on bounded domains in Hilbert spaces and extends them to finite-dimensional normed spaces with a strongly convex squared norm. For arbitrary nonempty closed convex domains in these spaces, a known gradient bound yields parameter-free regret guarantees. Surprisingly, on an unconstrained Hilbert space, the construction recovers the \ac{KT} coin-betting algorithm exactly. We also give a parameter-free algorithm for unconstrained domains that requires no gradient bound.
\end{abstract}

\maketitle

\begin{acronym}[FTRL] 
\acro{FTRL}[FTRL]{Follow-the-Regularized-Leader}
\acro{CBA}[CBA]{Conic Blackwell Algorithm}
\acro{CBA+}[CBA$^{+}$]{Conic Blackwell Algorithm Plus}
\acro{OMD}[OMD]{Online Mirror Descent}
\acro{OCO}[OCO]{Online Convex Optimization}
\acro{RM}[RM]{Regret Matching}
\acro{RM+}[RM$^{+}$]{Regret Matching$^{+}$}
\acro{KT}[KT]{Krichevsky--Trofimov}
\acro{SGD}[SGD]{Stochastic Gradient Descent}
\end{acronym}
\acresetall

\section{Introduction}
First-order \ac{OCO}~\citep{Gordon99,Gordon99b,Zinkevich03,Cesa-BianchiL06,Shalev-Shwartz12,Hazan16,Orabona19} and stochastic optimization algorithms often suffer from suboptimal empirical and theoretical performance when used with untuned step sizes. To remedy this situation, several strategies have been proposed: self-confident algorithms~\citep{AuerCG02}, \ac{RM} and \ac{RM+} on the simplex~\citep{HartMC00,TammelinBJB15}, AdaGrad step sizes~\citep{DuchiHS10,McMahanS10,StreeterM10}, parameter-free\footnote{See \citet{Orabona19} for an exact definition of this term in online learning.} online algorithms, including coin-betting algorithms for unbounded domains~\citep[e.g.,][]{ChaudhuriYH09,ChernovV10,StreeterM12,Orabona14,LuoS15,KoolenVE15,vanErvenK16,OrabonaP16,OrabonaT17,KempkaKW19,OrabonaP21,ChenKLO22,ZhangCP22,CutkoskyDM23,CutkoskyM24}, parameter-free \ac{SGD}~\citep{CarmonH22,IvgiHC23,ChzhenGS23}, and \ac{CBA} and \ac{CBA+} on compact convex sets~\citep{GrandClementK21,grand2024solving}. Notably, very few of these approaches share a common design strategy, with specialized assumptions and analyses being the norm.

The objective of this note is to introduce a \emph{unifying framework} for learning-rate-free algorithms for \ac{OCO}.
The companion paper by \citet{Kwon26approachability} proposes a general family of parameter-free algorithms for Blackwell's approachability
and applies them to unconstrained regret minimization using two key ingredients:
\begin{itemize}
    \item {\em Conic lifting:} add a scalar coordinate to the original
      (possibly unbounded) domain and take the conic hull of the lifted domain.
    \item {\em Regret minimization over cones:} apply suitably chosen
      regret minimizers, such as \ac{FTRL} or \ac{OMD},
      to this conic hull.
\end{itemize}
\citet{Kwon26approachability} focuses on unconstrained regret
minimization and Euclidean geometries. We extend this approach to general
(possibly unbounded) closed convex domains, non-Euclidean geometries,
and algorithms related to coin-betting.
The resulting learning-rate-free no-regret algorithms recover and extend \ac{CBA}, \ac{CBA+}, and \ac{KT} coin-betting~\citep{OrabonaP16}. Moreover, the construction directly incorporates non-Euclidean norms without black-box reductions~\citep{CutkoskyO18}.
Our contributions cover three settings.
\begin{enumerate}
\item \emph{Bounded domains.} We recover \ac{CBA} and \ac{CBA+} and extend them to non-Euclidean norms. The domain can be any nonempty closed bounded convex set in a real Hilbert space, or in a finite-dimensional normed space whose squared norm is strongly convex. The algorithms require a bound on the domain radius, but no horizon or gradient bound. Theorem~\ref{thm:compact-euclidean-regret} gives the Hilbert-space guarantee, and Theorem~\ref{thm:compact-leading-constant} states the sharpness of its leading constant, which improves on the corresponding AdaGrad constant.
\item \emph{Known gradient bound and closed convex domains.} We develop \ac{FTRL}-based
  algorithms for any nonempty closed convex domain, including unbounded ones.
  In Hilbert spaces, each update reduces to projections onto scaled copies of the translated feasible set and a scalar equation. On an unconstrained domain, the algorithm coincides with \ac{KT} coin-betting.
\item \emph{Unknown gradient bound on the whole space.} We give an adaptive \ac{FTRL}-based construction on the whole space under the same assumptions on the norm. Its regret guarantee includes a correction for increases in the largest observed gradient norm. The algorithm requires neither a gradient bound nor a horizon.
\end{enumerate}
We remark that all our regret upper bounds immediately yield guarantees for offline and stochastic optimization through the well-known online-to-batch conversion~\citep{Cesa-BianchiCG04,Orabona19}.

\subsection*{Outline} The rest of this paper is organized as follows. Section~\ref{sec:notation} introduces the notation and preliminary definitions. Section~\ref{sec:conic} presents the conic lifting, and Section~\ref{sec:rel} reviews the literature on conic lifting. Section~\ref{sec:bounded-euclidean} treats bounded domains in real Hilbert spaces. Section~\ref{sec:bounded-non-euclidean} extends the analysis to strongly convex squared norms on bounded domains, including the \ac{CBA+} construction. Section~\ref{sec:unbounded} treats possibly unbounded closed convex domains with a known gradient bound, and Section~\ref{sec:energy-adaptation} treats the unconstrained case without a supplied gradient bound.

\section{Notation and preliminaries}
\label{sec:notation}

We write vectors in bold lowercase, matrices in bold uppercase, and named sets in calligraphic letters. In finite dimensions, $\indicator$ is the all-ones vector. We identify a real Hilbert space $\mathcal{H}$ with its continuous dual and write $\|\cdot\|_{2}$ for its Hilbert norm, which is Euclidean when $\mathcal{H}=\mathbb{R}^{d}$.

\subsection*{Dual pairing and norms}
For a real normed space $\mathcal{X}$, its continuous dual $\mathcal{X}^{\star}$ consists of bounded linear functionals. We write $\langle\bg,\bx\rangle\coloneqq\bg(\bx)$ for $\bg\in\mathcal{X}^{\star}$ and $\bx\in\mathcal{X}$; under the Hilbert-space identification, this is the inner product. The dual norm is
\[
 \|\bg\|_{\star}
 \coloneqq\sup_{\|\bx\|\leqslant1} \ \langle\bg,\bx\rangle.
\]
A seminorm satisfies the norm axioms but may vanish on nonzero vectors.

\subsection*{Convex functions}
For a convex function $f$ on a real normed space $\mathcal{X}$ and $\by\in\operatorname{dom}f$, its subdifferential $\partial f(\by)$ consists of the $\bxi\in\mathcal{X}^{\star}$ such that, for every $\bx\in\operatorname{dom}f$,
\[
 f(\bx)
 \geqslant f(\by)+\langle\bxi,\bx-\by\rangle\,.
\]
For $\bxi\in\partial f(\by)$, the generalized Bregman divergence is
\[
 B_{f}^{\bxi}(\bx;\by)
 \coloneqq f(\bx)-f(\by)-\langle\bxi,\bx-\by\rangle\,.
\]
We call $f$ the distance-generating function and write $B_{f}$ when $\bxi=\nabla f(\by)$.
For $\mu>0$, the function $f$ is $\mu$-strongly convex with respect to a seminorm $|\cdot|$ if, for all $\bx,\by$ in its convex domain and $\lambda\in[0,1]$,
\[
 f\left((1-\lambda)\bx+\lambda\by\right)
 \leqslant(1-\lambda)f(\bx)+\lambda f(\by)-\frac{\mu}{2}\lambda(1-\lambda)|\bx-\by|^{2}\,.
\]
In this case, $B_{f}^{\bxi}(\bx;\by)\geqslant\mu|\bx-\by|^{2}/2$ for every $\bxi\in\partial f(\by)$.

For a convex set $\mathcal{V}\subseteq\mathcal{X}$ and $\bx\in\mathcal{V}$, its outward normal cone $\mathcal{N}_{\mathcal{V}}(\bx)$ consists of the $\bnu\in\mathcal{X}^{\star}$ such that, for every $\by\in\mathcal{V}$,
\[
 \langle\bnu,\by-\bx\rangle
 \leqslant0\,.
\]

\subsection*{Regret}
In the \ac{OCO} setting with a nonempty convex feasible set $\mathcal{V}\subseteq\mathcal{X}$, on round $t$ the learner chooses $\bx_{t}\in\mathcal{V}$ before observing a subgradient $\bg_{t}\in\partial\ell_{t}(\bx_{t})$ of the convex loss $\ell_{t}$. For a fixed comparator $\bu\in\mathcal{V}$ and horizon $T$, linearized regret is
\begin{equation}
\label{eq:regret-definition}
 \Regret_{T}(\bu)
 \coloneqq\sum_{t=1}^{T}\langle\bg_{t},\bx_{t}-\bu\rangle,
 \qquad
 \sum_{t=1}^{T}\left(\ell_{t}(\bx_{t})-\ell_{t}(\bu)\right)
 \leqslant\Regret_{T}(\bu)\,.
\end{equation}
We use \emph{gradient} as shorthand for the observed subgradient $\bg_{t}$. Throughout the paper, regret means linearized regret; its bounds also apply to convex loss regret.

We use the following elementary results; proofs are recalled in Appendix~\ref{app:prf-preliminaries}.

\begin{lemma}[Subgradients of a squared norm]
\label{lem:squared-norm-subgradient}
Let $\phi(\by)\coloneqq\|\by\|^{2}/2$ on a real normed space. Every $\bg\in\partial\phi(\by)$ satisfies $\langle\bg,\by\rangle=\|\by\|^{2}$ and $\|\bg\|_{\star}=\|\by\|$.
\end{lemma}

\begin{lemma}[Strong convexity and dual norm on a product space]
\label{lem:product-convexity}
Let $\mathcal{X}$ and $\mathcal{Y}$ be real normed spaces and let $\mu,\nu>0$. Assume that $f_{\mathcal{X}}\colon\mathcal{X} \to \R$ is $\mu$-strongly convex with respect to $\|\cdot\|_{\mathcal{X}}$ and $f_{\mathcal{Y}}\colon\mathcal{Y} \to \R$ is $\nu$-strongly convex with respect to $\|\cdot\|_{\mathcal{Y}}$. Then, $f((\bx,\by)) \coloneqq f_{\mathcal{X}}(\bx)+f_{\mathcal{Y}}(\by)$ is $1$-strongly convex with respect to the product norm below, with the indicated dual norm:
\[
\|(\bx,\by)\|
\coloneqq\sqrt{\mu\|\bx\|_{\mathcal{X}}^{2}+\nu\|\by\|_{\mathcal{Y}}^{2}}
\quad\text{and}\quad
\|(\bg,\bz)\|_{\star}
=\sqrt{\frac{\|\bg\|_{\mathcal{X},\star}^{2}}{\mu}+\frac{\|\bz\|_{\mathcal{Y},\star}^{2}}{\nu}},
\]
where $(\bx,\by)\in\mathcal{X}\times\mathcal{Y}$ and $(\bg,\bz)\in\mathcal{X}^{\star}\times\mathcal{Y}^{\star}$, and $\|\cdot\|_{\mathcal{X},\star}$ and $\|\cdot\|_{\mathcal{Y},\star}$ are the respective dual norms.
\end{lemma}

\begin{lemma}[Growth at a constrained minimum]
\label{lem:minimum-growth}
Let $\mathcal{D}$ be a convex subset of a real normed space and let $f\coloneqq p+r$ be finite on $\mathcal{D}$,
where $p$ is $\mu$-strongly convex with respect to a seminorm $|\cdot|$, with $\mu>0$,
and $r$ is convex and continuous on an open convex set containing $\mathcal{D}$.
Let $\iota_{\mathcal{D}}$ be zero on $\mathcal{D}$ and $+\infty$ outside it.
If $\bx^{\star}\in\argmin_{\bx\in\mathcal{D}}\ f(\bx)$, then there exists $\bxi\in\partial r(\bx^{\star})$ such that $-\bxi\in\partial(p+\iota_{\mathcal{D}})(\bx^{\star})$.
Every such $\bxi$ and every $\by\in\mathcal{D}$ satisfy
\begin{equation}
\label{eq:minimum-growth}
 f(\by)-f(\bx^{\star})
 \geqslant\frac{\mu}{2}|\by-\bx^{\star}|^{2}+r(\by)-r(\bx^{\star})-\langle\bxi,\by-\bx^{\star}\rangle\,.
\end{equation}
\end{lemma}

\begin{lemma}[Logarithm bounds]
\label{lem:logarithm-bound}
For every real $x<1$,
\begin{equation}
\label{eq:logarithm-general-bound}
 -\ln(1-x)\leqslant x+\frac{x^{2}}{2(1-\max\{x,0\})}\,.
\end{equation}
In particular, for $|x|\leqslant\frac{1}{2}$,
\begin{equation}
\label{eq:logarithm-upper-bound}
-\ln(1-x)
\leqslant x+x^{2}\,.
\end{equation}
Also, for $x>0$,
\begin{equation}
\label{eq:logarithm-ratio-bounds}
\ln x
\geqslant1-\frac{1}{x}\,.
\end{equation}
\end{lemma}

\section{Conic lifting}\label{sec:conic}
Conic lifting reduces regret minimization over a nonempty closed bounded convex set to regret minimization over a cone. We first give the reduction, then discuss its relation to earlier constructions.
\subsection{Lifted cone and regret reduction}

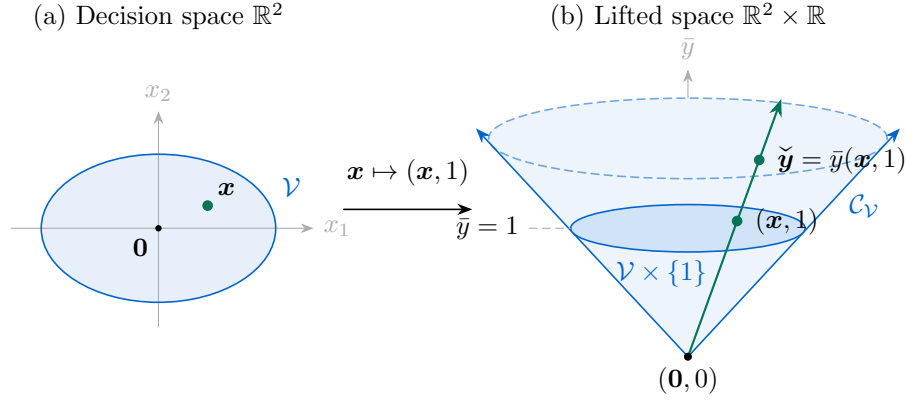
\begin{figure}[tbp]
\centering
\begin{tikzpicture}[font=\small, >=Stealth, line cap=round, line join=round]
 \colorlet{conicblue}{iroLink}
 \colorlet{conicaccent}{iroCite}
 \tikzset{
  conic axis/.style={->, thin, gray!65},
  conic boundary/.style={draw=conicblue, semithick},
  conic ray/.style={conicaccent, thick, ->}
 }
 \node at (0,4.5) {(a) Decision space $\R^{2}$};
 \node at (7,4.5) {(b) Lifted space $\R^{2}\times\R$};

 \begin{scope}[shift={(0,1.7)}]
  \fill[conicblue!10] (0,0) ellipse[x radius=1.55cm, y radius=0.98cm];
  \draw[conic axis] (-1.95,0) -- (2.05,0) node[right] {$x_{1}$};
  \draw[conic axis] (0,-1.3) -- (0,1.55) node[above] {$x_{2}$};
  \draw[conic boundary] (0,0) ellipse[x radius=1.55cm, y radius=0.98cm];
  \node[conicblue, anchor=west] at (1.5,0.55) {$\mathcal{V}$};
  \fill (0,0) circle[radius=1.2pt] node[below left] {$\bzero$};
  \fill[conicaccent] (0.65,0.3) circle[radius=2pt]
   node[above right, text=black] {$\bx$};
 \end{scope}

 \draw[->, semithick] (2.45,1.95) -- (4.15,1.95)
  node[midway, above=5pt] {$\bx\mapsto(\bx,1)$};

 \begin{scope}[shift={(7,0)}]
  \pgfmathsetmacro{\conicangle}{asin(-0.3136/1.7)}
  \pgfmathsetmacro{\conicotherangle}{180-\conicangle}
  \draw[conic axis] (0,0) -- (0,3.8) node[above] {$\bar{y}$};
  \path[fill=conicblue!7]
   (0,0) -- ({2.635*cos(\conicangle)},{2.89+0.53312*sin(\conicangle)})
   arc[start angle=\conicangle, end angle=\conicotherangle,
       x radius=2.635cm, y radius=0.53312cm] -- cycle;
  \draw[conic boundary, densely dashed, opacity=0.55]
   (0,2.89) ellipse[x radius=2.635cm, y radius=0.53312cm];
  \draw[conic boundary, ->] (0,0)
   -- ({1.84*1.55*cos(\conicangle)},{1.84*(1.7+0.3136*sin(\conicangle))});
  \draw[conic boundary, ->] (0,0)
   -- ({1.84*1.55*cos(\conicotherangle)},{1.84*(1.7+0.3136*sin(\conicotherangle))});
  \filldraw[fill=conicblue!18, conic boundary]
   (0,1.7) ellipse[x radius=1.55cm, y radius=0.3136cm];
  \draw[thin, gray!65, densely dashed] (-2.1,1.7) -- (-1.55,1.7);
  \node[anchor=east] at (-2.1,1.7) {$\bar{y}=1$};
  \node[conicblue] at (-0.35,1.1) {$\mathcal{V}\times\{1\}$};
  \node[conicblue, anchor=west] at (2,2) {$\mathcal{C}_{\mathcal{V}}$};
  \draw[conic ray] (0,0) -- (1.235,3.4124);
  \fill[conicaccent] (0.65,1.796) circle[radius=2pt]
   node[right=3pt, text=black] {$(\bx,1)$};
  \fill[conicaccent] (0.9425,2.6042) circle[radius=2pt]
   node[right=3pt, text=black] {$\widecheck{\by}=\bar{y}(\bx,1)$};
  \fill (0,0) circle[radius=1.5pt] node[below] {$(\bzero,0)$};
 \end{scope}
\end{tikzpicture}
\caption{Conic lifting of a bounded closed convex set $\mathcal{V}\subseteq\R^{2}$.
The slice at $\bar{y}=1$ is $\mathcal{V}\times\{1\}$.
Every positive multiple $\widecheck{\by}=\bar{y}(\bx,1)$ on the highlighted ray decodes to the same decision $\bx=\by/\bar{y}$.
The cone continues along the indicated rays beyond the displayed portion.}
\label{fig:conic-lifting}
\end{figure}

Let $\mathcal{V}\subseteq\mathcal{H}$ be nonempty, closed, bounded, and convex, where $\mathcal{H}$ is a real Hilbert space.
Fix $R>0$ such that $\|\bx\|_{2}\leqslant R$ for every $\bx\in\mathcal{V}$, and choose $\bx_{0}\in\mathcal{V}$. We define the cone
\begin{equation}
\label{eq:compact-cone}
 \mathcal{C}_{\mathcal{V}}
 \coloneqq\{\lambda(\bx,1) : \bx\in\mathcal{V},\ \lambda\geqslant0 \}\,.
\end{equation}
A check mark denotes a lifted vector: $\widecheck{\by}=(\by,\bar{y})\in\mathcal{H}\times\mathbb{R}$. Lifted losses and subgradients belong to the product dual $\mathcal{H}^{\star}\times\mathbb{R}$.
Since $\mathcal{V}$ is closed and bounded, $\cCV$ is a nonempty closed convex cone, so orthogonal projection onto $\cCV$ is well defined in $\mathcal{H}\times\R$. The cone $\cCV$ is sometimes called the \emph{homogenization cone} of $\cV$~\citep[see, e.g., Section~1.5 in][]{brinkhuis2020convex}, in part because its orthogonal projection is homogeneous of degree $1$: $\Pi_{\cCV}(\eta \widecheck{\by}) = \eta \Pi_{\cCV}(\widecheck{\by})$ for every $\widecheck{\by} \in \mathcal{H} \times \mathbb{R}$ and $\eta>0$.
Figure~\ref{fig:conic-lifting} illustrates the construction for a two-dimensional set.

For $\widecheck{\by}=(\by,\bar{y})\in\mathcal{C}_{\mathcal{V}}$, define its decoded decision by
\begin{equation}\label{eq:decode}
 \decode(\widecheck{\by})\coloneqq
 \begin{cases}
  \by/\bar{y},&\text{if $\bar{y}>0$},\\
  \bx_{0},&\text{otherwise},
 \end{cases}
\end{equation}
where $\bx_{0}\in\mathcal{V}$ is the fixed decision chosen above.
Given a decision $\bx\in\mathcal{V}$ and a loss vector $\bg\in\mathcal{H}^{\star}$, define the lifted loss vector in the product dual $\mathcal{H}^{\star}\times\mathbb{R}$ by
\begin{equation}\label{eq:lifted loss}
 \widecheck{\bg}
 \coloneqq(\bg,-\langle\bg,\bx\rangle)\,.
\end{equation}
A crucial property is that, for every $\widecheck{\by}\in\mathcal{C}_{\mathcal{V}}$ and $\bg\in\mathcal{H}^{\star}$, the lifted loss associated with $\bx=\decode(\widecheck{\by})$ satisfies the orthogonality identity
\begin{equation}\label{eq:orthogonality blackwell}
 \langle\widecheck{\bg},\widecheck{\by}\rangle
 =0\,.
\end{equation}
Indeed, if $\bar{y}>0$, then
\[
 \langle\widecheck{\bg},\widecheck{\by}\rangle
 =\langle\bg,\by\rangle-\bar{y}\langle\bg,\by/\bar{y}\rangle
 =0\,.
\]
If $\bar{y}=0$, then $\widecheck{\by}=(\bzero,0)$, so the same identity holds.

The following simple lemma illustrates the key idea behind conic lifting: a reduction from regret minimization over $\mathcal{V}$ to regret minimization over the lifted cone $\mathcal{C}_{\mathcal{V}}$.
\begin{lemma}\label{lem: regret reduction}
Let $\widecheck{\by}_{1},\dots,\widecheck{\by}_{T}\in\mathcal{C}_{\mathcal{V}}$ and $\bg_{1},\dots,\bg_{T}\in\mathcal{H}^{\star}$.
For each $t$, set $\bx_{t}\coloneqq\decode(\widecheck{\by}_{t})$ and $\widecheck{\bg}_{t}\coloneqq(\bg_{t},-\langle\bg_{t},\bx_{t}\rangle)$.
Then, for every $\bu\in\mathcal{V}$, with $\widecheck{\bu}\coloneqq(\bu,1)$,
\[
 \sum_{t=1}^{T}\langle\bg_{t},\bx_{t}-\bu\rangle
 =\sum_{t=1}^{T}\langle\widecheck{\bg}_{t},\widecheck{\by}_{t}-\widecheck{\bu}\rangle\,.
\]
\end{lemma}
The same reduction appears in Section~3 of \citet{gordon2006no}, between equations~(2) and~(3), Lemma~2.1 of \citet{farina2023regret}, and Proposition~3.1 of \citet{ChakrabartiGK24}.
\begin{proof}[Proof of Lemma \ref{lem: regret reduction}]
By the definition of $\widecheck{\bg}_{t}$ and the orthogonality identity \eqref{eq:orthogonality blackwell}, respectively,
\[
 \sum_{t=1}^{T}\langle\bg_{t},\bx_{t}-\bu\rangle
 =\sum_{t=1}^{T}\langle-\widecheck{\bg}_{t},\widecheck{\bu}\rangle
 =\sum_{t=1}^{T}\langle\widecheck{\bg}_{t},\widecheck{\by}_{t}-\widecheck{\bu}\rangle\,.
\]
\end{proof}

\section{Related work}
\label{sec:rel}

\subsection*{Regret via approachability}
Lemma~\ref{lem: regret reduction} follows from the definition of the
lifted loss in \eqref{eq:lifted loss} and the orthogonality property
\eqref{eq:orthogonality blackwell}. This orthogonality can be interpreted as a hyperplane forcing property in Blackwell
approachability~\citep{Blackwell56}. With uniformly bounded gradients,
the bounded-domain regret minimizers below ensure that the average negative lifted loss vector
$-\frac{1}{T}\sum_{t=1}^{T} \widecheck{\bg}_{t}$ approaches the polar
cone of $\cCV$ as $T \to + \infty$. The use of Blackwell's approachability to minimize regret \emph{over the simplex} by adding a
dimension was noted very early in~\citet{blackwell1954controlled}; see
\citet{perchet2014approachability} for the extension to a convex
compact domain.
\citet{Kwon26approachability} applies conic lifting to unconstrained
regret minimization by constructing regularizers whose convex
conjugates have gradients that avoid the zero-height boundary of the
closure of the lifted cone, except at the origin. The present paper
extends this approach to general unbounded closed convex domains and
non-Euclidean squared norms in Sections~\ref{sec:unbounded} and
\ref{sec:energy-adaptation}.

\subsection*{Approachability via regret minimization on cones}
The present work relies on regret minimizers on the lifted cone, and the
lineage of this idea will be discussed in paragraphs below. But
given the above approachability interpretation, this
approach is a special case of converting regret
minimizers into approachability algorithms for conic targets. This 
general approach first appeared in~\citet{abernethy2011blackwell}
which focuses on minimizing the Euclidean distance to the target by
defining the regret minimizers on the intersection of the relevant
closed convex cone with the unit Euclidean ball. This conversion
scheme was refined in \citet{kwon2016mirror,kwon2021refined} which
allowed for bounds on more general non-Euclidean distance-like
quantities by defining regret minimizers on any convex compact subset
that contains all directions of the cone. In the approach presented by \citet{Kwon26approachability} regret minimizers are
instead defined on the whole cone (which, interestingly, is no longer
required to be closed) with the help of near-2-positively-homogeneous
regularizers which yield parameter-free guarantees.

\subsection*{Simplex domains}
The original regret minimizer of \citet{blackwell1954controlled}
uses approachability for a lifted set slightly different from $\cCV$
(see \citet[p.\ 337]{blackwell1954controlled}).
For the simplex, however, approachability can be applied without lifting:
one can use the auxiliary loss vectors
$\bg_{t} - \langle \bg_{t},\bx_{t}\rangle \indicator$ instead of
$\widecheck{\bg}_{t}$, together with $\R_{+}^{d}$, the
conic hull of the simplex $\Delta^{d}$, without
introducing an additional coordinate~\citep{HartMC00}. Running \ac{FTRL}
with a Euclidean regularizer on $\R_{+}^{d}$ to minimize regret over the
simplex in $\R^{d}$ yields \ac{RM}~\citep{HartMC00}, as noted by
Nicol\`{o} Cesa-Bianchi in an unpublished manuscript titled ``The Joys of
Regret Matching'', dated October 14, 2015. More generally, the equivalence between Euclidean \ac{FTRL} for conic targets and Blackwell's algorithm was independently established in \citet{kwon2016mirror,kwon2021refined}. To improve practical performance in learning in games,
\citet{Tammelin2014,TammelinBJB15} proposed \ac{RM+}, a variant using
thresholding. \citet{FarinaKS21} and, independently,
\citet{FlaspohlerOCMOOM21} identified it as \ac{OMD} with a Euclidean
regularizer on $\R_{+}^{d}$. In particular, the inverse $\ell_{1}$-norm of the nonnegative
regret vector maintained by \ac{RM+} can be interpreted as an effective
step size~\citep[Section~4.1]{zhang2025scale}.

\subsection*{Lifting general bounded convex domains}
For a compact convex set $\cV$, the same conic techniques can be applied
directly in the ambient space $\mathcal{H}$, without adding a dimension,
whenever there exists $\be\in\mathcal{H}^{\star}$ such that
$\langle\be,\bx\rangle=1$ for every $\bx\in\cV$. One then uses the conic
hull of $\cV$ instead of $\cCV$, with loss vectors
$\bg-\langle\bg,\bx\rangle\be\in\mathcal{H}^{\star}$. The simplex is the
special case $\be=\indicator$. \citet{ChakrabartiGK24} use this construction
for a sequence-form polytope $\cV$ in learning in extensive-form games.

However, when no such normalizing vector is available, lifting
$\bx \mapsto (\bx,1)$ from $\mathcal{H} $ to $\mathcal{H} \times \R$
provides the \emph{lifted copy} $\{(\bx,1) \; | \; \bx \in \cV\}$ to which the same
techniques apply. The normalizing vector is then
$(\boldsymbol{0},1)\in\mathcal{H}^{\star}\times\R$, since
$\{(\bx,1) \; | \; \bx \in \cV\} \subset \{\bz \in \mathcal{H} \times \R \; | \; \langle \bz,(\boldsymbol{0},1) \rangle = 1\}$.
This is the main point of using the additional coordinate and the lifted cone $\cCV$.

The reduction from regret minimization over $\cV$ to regret minimization
over $\cCV$ has been rediscovered several times.
\citet{gordon2005no,gordon2006no} proposed {\em Lagrangian Hedging} for a
compact convex decision set $\cV$. This framework can be interpreted as
running \ac{FTRL} on $\cCV$, with particular attention to Euclidean and
entropic regularizers.
Note, however, that the lifting is implicit in \citet[Section 3]{gordon2006no}, which works
in a finite-dimensional space and assumes the existence of a vector $\be$ such that
$\langle\bx,\be\rangle=1$ for all $\bx\in\cV$. Gordon briefly explains that
$\cCV$ can be used instead when no such vector exists.
The same conic lifting and regret reduction appear in
\citet[Theorem~16]{abernethy2011blackwell}, who reduce regret
minimization to approachability. Works that depart from \ac{FTRL}
include \citet{GrandClementK21} (see also the journal
version by \citet{grand2024solving}) and \citet{d2021optimistic}.
\citet{GrandClementK21} introduce the \ac{CBA}, which runs \ac{FTRL}
with a Euclidean regularizer on the cone $\cCV$ to minimize regret on
a convex compact set $\cV$, and the \ac{CBA+}, which runs \ac{OMD}
with a Euclidean regularizer on $\cCV$ to minimize regret on $\cV$.
Their analysis focuses on Blackwell approachability; Appendix~A of
\citet{GrandClementK21} discusses the connection to regret minimization
on the cone. For the whole-cone \ac{FTRL}/\ac{OMD} viewpoint, see also
\citet[Sections~4--5]{FarinaKS21}. \citet{d2021optimistic} generalize
Lagrangian Hedging to adaptive step sizes, optimistic updates, and
internal/$\Phi$-regret. For the simplex, they introduce Lagrangian
Hedging$^{+}$ (see their page 6), a version of \ac{RM+} for $\Phi$-regret:
\ac{OMD} with a Euclidean regularizer on $\R_{+}^{d}$. As noted above, no additional coordinate is needed
before taking the conic hull.

\subsection*{Truncated lifts and learning in games}
Conic constructions also yield uncoupled algorithms with near-constant ($O(\ln T)$) per-player external regret in multiplayer games with finite action sets (see Section~4 of \citet{soleymani2025cautious}), and with constant per-player external regret (see Section~3.1 of \citet{abbadi2026constant} and Section~3.1 of \citet{liu2026constant}). Interestingly, on the simplex, COFTRL~\citep{soleymani2025cautious} and HOOD~\citep{abbadi2026constant} use \emph{truncated} conic hulls without adding a coordinate. ECHO-OFTRL~\citep{liu2026constant} instead uses the added-coordinate lift $\{\lambda (\bx,1) \; | \; \bx \in \Delta^{d},\lambda \in [0,1]\}$; COFTRL also uses an added coordinate in its extension to general compact convex domains. The closed comparator sets for these truncated constructions include zero, whereas $\cCV$ allows every $\lambda \geqslant 0$. Section~\ref{sec:energy-adaptation} uses a different truncated cone. 

Truncated conic lifting for games originates in \citet{farina2022near}, who obtain near-constant regrets for general convex games using optimistic \ac{FTRL} with a logarithmic regularizer on the truncated lifted cone. Here, the lifting is useful because the lifted regrets are nonnegative, whereas large negative original regrets for some players can counterbalance large positive regrets for others, complicating the analysis. \citet{farina2022near} trace truncated conic lifting to \citet{lee2020bias,luo2022adaptive} for online learning with bandit feedback. In turn, \citet{lee2020bias} traces lifting to \citet{nesterov1994interior}, whose Proposition~5.1.4 constructs normal barriers for the lifted cone from a self-concordant barrier for the domain.

Interestingly, a third recent algorithm with constant per-player regret in finite multiplayer games uses Blackwell approachability through \ac{RM}-style regret minimizers on simplexes~\citep{soleymani2026multiplicative}.

\begin{table}[htbp]
\centering
\begin{threeparttable}
\caption{Uses of conic reformulations in online learning.
\emph{Lifting} means appending a coordinate, $x\mapsto(x,1)$. \ac{RM} and \ac{RM+} are from
\citet{HartMC00,TammelinBJB15}. HOOD is from \citet{abbadi2026constant}. ECHO-OFTRL is from 
\citet{liu2026constant}. Lagrangian Hedging is from \citet{gordon2005no,gordon2006no}, and \ac{CBA}/\ac{CBA+} are from 
\citet{GrandClementK21,grand2024solving}. LRL-OFTRL is from \citet{farina2022near}, and COFTRL is from \citet{soleymani2025cautious}, who also treat compact convex domains using the lifted cone $\cCV$. For truncated domains, the regret identities use the closed comparator set, including zero.
}
\label{tab:conic-lifting-literature}
\small
\setlength{\tabcolsep}{4pt}
\renewcommand{\arraystretch}{1.12}
\begin{tabularx}{\linewidth}{@{}
  >{\raggedright\arraybackslash}p{0.15\linewidth}
  >{\raggedright\arraybackslash}p{0.15\linewidth}
  >{\centering\arraybackslash}p{0.085\linewidth}
  >{\raggedright\arraybackslash}p{0.30\linewidth}
  >{\raggedright\arraybackslash}X
  }
\toprule
\textbf{Papers / algorithms}
& \textbf{Decision set $\mathcal{V}$}
& \textbf{Lifting}
& \textbf{Cone (possibly truncated)}
& \textbf{Role of the conic construction} \\
\midrule

\textbf{\ac{RM} / \ac{RM+}}
& Simplex
& $\times$
& 
  \(\{\lambda x:x\in\Delta^{d},\ \lambda\geqslant0\}=\mathbb{R}_{+}^{d}\)
& Homogeneity and normalization yield learning-rate-free regret bounds. \\
\addlinespace[0.45em]

{\bf COFTRL}
& Simplex
& $\times$
&  \(\{\lambda x:x\in\Delta^{d},\ 0<\lambda\leqslant1\}\) 
& Lifted regret = positive part of maximum linearized regret. \\
\addlinespace[0.45em]

{\bf HOOD}
& Simplex
& $\times$
& 
  \(\{\lambda x:x\in\Delta^{d},\ 0\leqslant\lambda\leqslant1\}\) 
& Lifted regret = positive part of maximum linearized regret. \\
\addlinespace[0.45em]

{\bf ECHO-OFTRL}
& Simplex
& $\checkmark$
& 
  \(\{\lambda(x,1):x\in\Delta^{d},\ 0\leqslant\lambda\leqslant1\}\) 
& Lifted regret = positive part of maximum linearized regret. \\

\midrule
{\bf LRL-OFTRL}
& Compact convex
& $\checkmark$
& \(\{\lambda(x,1):x\in\mathcal{V},\ 0\leqslant\lambda\leqslant1\}\)
& Lifted regret = positive part of maximum linearized regret.\\

\citet{lee2020bias,luo2022adaptive}
& Compact convex
& $\checkmark$
&   \(\{\lambda(x,1):x\in\mathcal{V},\ \lambda\geqslant0\}\)
& A logarithmically homogeneous barrier controls stability as learning rates change. \\

\textbf{Lagrangian Hedging,}
\par \textbf{\ac{CBA} / \ac{CBA+}}
& Compact convex
& $\checkmark$
& 
  \(\{\lambda(x,1):x\in\mathcal{V},\ \lambda\geqslant0\}\)
  & Orthogonality and rescaling lifted comparators yield learning-rate-free regret bounds.
\\

\midrule

  \citet{Kwon26approachability} and
\textbf{this paper}\par
& Closed convex
& $\checkmark$
& \(\{\lambda(x,1):x\in\mathcal{V},\ \lambda\geqslant0\}\)
& Orthogonality and rescaling lifted comparators yield learning-rate-free regret bounds.
\\
\bottomrule
\end{tabularx}
\end{threeparttable}
\end{table}

\section{Warm-up: the Conic Blackwell Algorithm}\label{sec:bounded-euclidean}
We begin with the \ac{CBA} of \citet{GrandClementK21}, originally formulated over compact convex sets in Euclidean space. Lemma~\ref{lem: regret reduction} reduces the construction to regret minimization over $\mathcal{C}_{\mathcal{V}}$, where \ac{CBA} uses \ac{FTRL}. We give the bounded Hilbert-space guarantee; Section~\ref{sec:bounded-non-euclidean} treats strongly convex squared norms and the \ac{OMD} formulation of \ac{CBA+}. To account for the radius $R$, we equip the lifted Hilbert space with the norm and dual norm
\begin{equation}
\label{eq:compact-metric}
\|\widecheck{\by}\|_{R}^{2}
= \|(\by,\bar{y})\|_{R}^{2}
\coloneqq\|\by\|_{2}^{2}+R^{2}\bar{y}^{2}
\quad\text{and}\quad
\|\widecheck{\bg}\|_{R,\star}^{2}
= \|(\bg,\bar{g})\|_{R,\star}^{2}
\coloneqq\|\bg\|_{2}^{2}+\frac{\bar{g}^{2}}{R^{2}}\,.
\end{equation}

We use \ac{FTRL} with a fixed quadratic regularizer on this cone.

\begin{algorithm}[t]
\caption{Conic Blackwell Algorithm}
\label{alg:compact-ftrl}
\algrenewcommand{\algorithmicrequire}{\textbf{Input:}}
\begin{algorithmic}[1]
\Require A feasible set $\mathcal{V}\subseteq\mathcal{X}$, a radius $R>0$, and $\bx_{0}\in\mathcal{V}$
\State $\mathcal{C}_{\mathcal{V}}\coloneqq\{\lambda(\bx,1):\bx\in\mathcal{V},\ \lambda\geqslant0\}$
\State $\widecheck{\by}_{1}=(\by_{1},\bar{y}_{1})\coloneqq(\bzero,0)$
\For{$t=1,2,\dots$}
 \State $\bx_{t} = \decode(\displaystyle\widecheck{\by}_{t})$
 \State Play $\bx_{t}$ and observe $\bg_{t}\in\partial\ell_{t}(\bx_{t})$
 \State $\widecheck{\bg}_{t}\coloneqq(\bg_{t},-\langle\bg_{t},\bx_{t}\rangle)$
 \State $\displaystyle\widecheck{\by}_{t+1}=(\by_{t+1},\bar{y}_{t+1})\coloneqq\argmin_{\widecheck{\by}=(\by,\bar{y})\in\mathcal{C}_{\mathcal{V}}}\left\{\frac{1}{2}\left(\|\by\|_{\mathcal{X}}^{2}+R^{2}\bar{y}^{2}\right)+\sum_{i=1}^{t}\langle\widecheck{\bg}_{i},\widecheck{\by}\rangle\right\}$ \label{step: algo 1 FTRL}
\EndFor
\end{algorithmic}
\end{algorithm}

For Algorithms~\ref{alg:compact-ftrl} and~\ref{alg:compact-md}, let $\mathcal{X}$ be a finite-dimensional real normed space, or a real Hilbert space with its Hilbert norm, such that $\|\cdot\|_{\mathcal{X}}^{2}/2$ is $\mu$-strongly convex with respect to $\|\cdot\|_{\mathcal{X}}$ for some $\mu>0$. Let $\mathcal{V}\subseteq\mathcal{X}$ be nonempty, closed, bounded, and convex, with $\|\bx\|_{\mathcal{X}}\leqslant R$ for every $\bx\in\mathcal{V}$.
In Algorithm~\ref{alg:compact-ftrl}, take $\mathcal{X}=\mathcal{H}$ with its Hilbert norm. With the fixed lifted regularizer $\psi_{t}(\widecheck{\by})=\|\widecheck{\by}\|_{R}^{2}/2$, define
\[
 F_{t}(\widecheck{\by})\coloneqq\psi_{t}(\widecheck{\by})+\sum_{i=1}^{t-1}\langle\widecheck{\bg}_{i},\widecheck{\by}\rangle\,.
\]
The update minimizes $F_{t}$ over $\mathcal{C}_{\mathcal{V}}$; its feedback satisfies \eqref{eq:orthogonality blackwell}.

To bound regret on the cone, we use the classic \ac{FTRL} regret equality in Theorem~\ref{thm:ftrl-equality}.
The initial condition is $\psi_{T+1}(\widecheck{\bu})$, since $\psi_{1}(\widecheck{\by}_{1})=0$.
Since $F_{t}$ is $1$-strongly convex and minimized at $\widecheck{\by}_{t}$, Lemma~\ref{lem:minimum-growth} with $p=F_{t}$, $r=0$, and the norm $\|\cdot\|_{R}$ gives
\[
 F_{t}\left(\widecheck{\by}_{t}\right)-F_{t+1}\left(\widecheck{\by}_{t+1}\right)+\left\langle\widecheck{\bg}_{t},\widecheck{\by}_{t}\right\rangle
 \leqslant\frac{1}{2}\|\widecheck{\bg}_{t}\|_{R,\star}^{2}\,.
\]
The final condition $F_{T+1}\left(\widecheck{\by}_{T+1}\right)-F_{T+1}\left(\widecheck{\bu}\right)$ is nonpositive for every feasible $\widecheck{\bu}$. Substitution into the equality gives, for every $\widecheck{\bu}\in\mathcal{C}_{\mathcal{V}}$,
\[
\sum_{t=1}^{T}\langle\widecheck{\bg}_{t},\widecheck{\by}_{t}-\widecheck{\bu}\rangle
\leqslant\frac{1}{2}\|\widecheck{\bu}\|_{R}^{2} + \frac{1}{2}\sum_{t=1}^{T}\|\widecheck{\bg}_{t}\|_{R,\star}^{2}\,.
\]
To recover a regret bound on the original linear losses, for $\bu\in\mathcal{V}$, define $\widecheck{\bu}\coloneqq(\bu,1)$ and compare with $\lambda\widecheck{\bu}$ for $\lambda>0$.
From $\langle \widecheck{\bg}_{t}, \widecheck{\by}_{t}\rangle=0$, the left-hand side becomes $\lambda\Regret_{T}\left(\bu\right)$.
Following the comparator-rescaling argument of \citet{Kwon26approachability}, we divide by $\lambda$ and optimize to obtain
\[
 \Regret_{T}\left(\bu\right)\leqslant\inf_{\lambda>0}
 \left\{\frac{\lambda}{2}\|\widecheck{\bu}\|_{R}^{2}+\frac{1}{2\lambda}\sum_{t=1}^{T}\|\widecheck{\bg}_{t}\|_{R,\star}^{2}\right\}
 =\|\widecheck{\bu}\|_{R}\sqrt{\sum_{t=1}^{T}\|\widecheck{\bg}_{t}\|_{R,\star}^{2}}\,.
\]
Since
\[
 \|\widecheck{\bu}\|_{R}^{2}
 =\|\bu\|_{2}^{2}+R^{2}
 \leqslant2R^{2}\quad
 \text{and}\quad
 \|\widecheck{\bg}_{t}\|_{R,\star}^{2}
 =\|\bg_{t}\|_{2}^{2}+\frac{\langle\bg_{t},\bx_{t}\rangle^{2}}{R^{2}}
 \leqslant2\|\bg_{t}\|_{2}^{2},
\]
we obtain the following guarantee.
\begin{theorem}[Bounded-domain regret in a Hilbert space]
\label{thm:compact-euclidean-regret}
Run Algorithm~\ref{alg:compact-ftrl} in a real Hilbert space $\mathcal{X}\coloneqq\mathcal{H}$ with its Hilbert norm. Then, for every $T\geqslant1$ and every $\bu\in\mathcal{V}$,
\begin{equation}
\label{eq:compact-euclidean-regret}
\Regret_{T}(\bu)
\leqslant 2R\sqrt{\sum_{t=1}^{T}\|\bg_{t}\|_{2}^{2}}\,.
\end{equation}
\end{theorem}

\begin{remark}
Notice that initially, $\bar{y}_{1}=0$, and the condition $\bar{y}_{t}=0$ can persist for many rounds. However, once the iterate leaves the origin, it cannot return. Since $F_{t}(\bzero)=0$ and the minimizer is unique, $F_{t}(\widecheck{\by}_{t})<0$ if and only if $\widecheck{\by}_{t}\neq\bzero$. Moreover,
\[
F_{t+1}(\widecheck{\by}_{t+1})
=\min_{\widecheck{\by}\in\mathcal{C}_{\mathcal{V}}}\ F_{t+1}(\widecheck{\by})
\leqslant F_{t+1}(\widecheck{\by}_{t})
=F_{t}(\widecheck{\by}_{t})+\langle \widecheck{\bg}_{t},\widecheck{\by}_{t}\rangle
=F_{t}(\widecheck{\by}_{t})\,.
\]
Hence, once $\bar{y}_{t}\neq0$, the next iterates $\bar{y}_{t'}$ for $t'>t$ cannot be zero.
\end{remark}

\subsection{Comparison with AdaGrad-Norm}
AdaGrad-Norm~\citep{StreeterM10} instead gives, for every $\bu\in\mathcal{V}$,
\[
\Regret_{T}(\bu)
\leqslant 2 R \sqrt{2 \sum_{t=1}^{T}\|\bg_{t}\|_{2}^{2}}\,.
\]
Hence, Theorem~\ref{thm:compact-euclidean-regret} improves this leading constant to $2$. The following theorem shows that this constant is optimal for a bound of this form.

\begin{theorem}
\label{thm:compact-leading-constant}
Let $x_{t}$, $g_{t}$, and $u$ denote scalar predictions, gradients, and comparators. Fix $R>0$ and $0\leqslant c<2$. On $\mathcal{V}=[-R, R]$, no online algorithm, given neither a gradient bound nor a horizon, can guarantee, for every $T\geqslant1$, $u\in\mathcal{V}$, and gradient sequence,
\[
 \E\left[\Regret_{T}(u)\right]\leqslant cR\sqrt{\sum_{t=1}^{T}g_{t}^{2}}\,.
\]
The expectation is over any internal randomization of the algorithm.
The impossibility holds even when every gradient sequence is bounded by some finite, unknown constant.
\end{theorem}

\begin{proof}
Suppose such a guarantee holds. Fix any deterministic gradient history through round $t-1$. For any $M>0$, extend it by $g_{t}=-M$ and then zero gradients. The decision $x_{t}$ is chosen before observing $g_{t}$, so its distribution is independent of $M$. Against $u=R$, the assumed bound gives
\[
 \sum_{i=1}^{t-1}g_{i}\left(\E[x_{i}]-R\right)+M\left(R-\E[x_{t}]\right)
 \leqslant cR\sqrt{\sum_{i=1}^{t-1}g_{i}^{2}+M^{2}}\,.
\]
We divide by $M$ and take the limit as $M\to\infty$. The terms from the fixed history tend to zero, yielding $\E[x_{t}]\geqslant(1-c)R$ after every history.
Each extension used here has a finite-gradient bound; no single unbounded sequence is needed.

Now take $g_{t}=1$ on every round and compare with $u=-R$. The necessary restriction on every prediction implies
\[
 \E\left[\Regret_{T}(-R)\right]
 =\sum_{t=1}^{T}\left(\E[x_{t}]+R\right)
 \geqslant (2-c)RT\,.
\]
For $T>(c/(2-c))^{2}$, this exceeds the claimed bound $cR\sqrt{T}$, a contradiction.
\end{proof}

\subsection{Computing the update} We explain how to compute the update in Step~\ref{step: algo 1 FTRL} of Algorithm~\ref{alg:compact-ftrl}.
Let $\Pi_{\mathcal{C}_{\mathcal{V}}}^{R}$ denote projection in the norm $\|\cdot\|_{R}$. Completing the square in the update gives
\[
 \widecheck{\by}_{t}
 =\Pi_{\mathcal{C}_{\mathcal{V}}}^{R}
 \!\left(\left(-\sum_{i=1}^{t-1} \bg_{i},\frac{\sum_{i=1}^{t-1} \langle\bg_{i}, \bx_{i}\rangle}{R^{2}}\right)\right)\,.
\]
Multiplying $\bar{y}$ by $R$ turns this into an ordinary Hilbert-space projection onto the cone
\[
 \left\{\lambda\left(\bx,R\right):\lambda\geqslant0,\ \bx\in\mathcal{V}\right\}\,.
\]
Divide the last coordinate of the projection by $R$ afterward. The map $(\by,\bar{y})\mapsto(\by,-\bar{y})$ is an isometry, so the projection bounds cited below also apply to our convention $\bar{y}\geqslant0$.

In finite-dimensional Euclidean spaces, Moreau's decomposition reduces these weighted projections to projections onto the polar cone. For any closed convex cone $\mathcal{C}$, let $\Pi_{\mathcal{C}}$ denote the Euclidean projection. Moreau's decomposition theorem gives $\Pi_{\mathcal{C}}\left(\widecheck{\bv}\right) + \Pi_{\mathcal{C}^{\circ}}\left(\widecheck{\bv}\right) = \widecheck{\bv}$, where $\mathcal{C}^{\circ}\coloneqq\left\{\widecheck{\bz}:\left\langle \widecheck{\bz},\widecheck{\bx}\right\rangle\leqslant0,\ \forall \widecheck{\bx}\in\mathcal{C}\right\}$ is the polar cone. Hence, $\Pi_{\mathcal{C}}\left(\widecheck{\bv}\right)$ can be recovered from $\Pi_{\mathcal{C}^{\circ}}\left(\widecheck{\bv}\right)$. See \citet{GrandClementK21} for projection examples (the simplex, an $\ell_{p}$-ball with $p \in \{1,2,+\infty\}$, and the intersection of an $\ell_{2}$-ball with the simplex).

\section{Bounded domains with strongly convex squared norms}\label{sec:bounded-non-euclidean}
We now run the bounded-domain construction in a finite-dimensional normed space whose squared norm is strongly convex.

Let $\mathcal{X}$ be a finite-dimensional real normed space, and let $\mathcal{V}\subseteq\mathcal{X}$ be nonempty, convex, and compact. The lifting takes place in the product space $\mathcal{X}\times\R$.
Let $R>0$ satisfy $R\geqslant \|\bx\|_{\mathcal{X}}$ for all $\bx \in \mathcal{V}$. Assume $\psi_{\mathcal{X}}\colon\mathcal{X}\to \R$ defined by $\psi_{\mathcal{X}}(\bx)\coloneqq\frac{1}{2} \|\bx\|^{2}_{\mathcal{X}}$ is $\mu$-strongly convex with respect to $\|\cdot\|_{\mathcal{X}}$, with $\mu>0$.
Define the lifted space norm and its dual by
\begin{equation}
\label{eq:compact-product-norm}
\|\widecheck{\by}\|^{2}
= \|(\by,\bar{y})\|^{2}
\coloneqq\mu\|\by\|_{\mathcal{X}}^{2}+R^{2}\bar{y}^{2}
\quad\text{and}\quad
\|\widecheck{\bg}\|_{\star}^{2}
=\|(\bg,\bar{g})\|_{\star}^{2}
\coloneqq\frac{1}{\mu}\|\bg\|_{\mathcal{X},\star}^{2}+\frac{\bar{g}^{2}}{R^{2}}\,.
\end{equation}
Define the regularizer $\psi\colon\mathcal{X}\times\R \to \R$ as
\begin{equation}
 \label{eq:compact-product-potential}
 \psi(\widecheck{\by})
 = \psi((\by,\bar{y}))
 \coloneqq\frac{1}{2}\left(\|\by\|_{\mathcal{X}}^{2}+R^{2}\bar{y}^{2}\right)\,.
\end{equation}
For example, on $\mathbb{R}^{d}$ with $d\geqslant2$ and $1<p<2$, the function $\|\cdot\|_{p}^{2}/2$ is $(p-1)$-strongly convex with respect to $\|\cdot\|_{p}$~\citep[Lemma~6.35]{Orabona19}. This gives a non-Hilbert norm satisfying the assumption.

The two terms have strong-convexity constants $\mu$ and $R^{2}$ in their respective coordinates. Thus, Lemma~\ref{lem:product-convexity} gives the dual norm in \eqref{eq:compact-product-norm} and makes $\psi$ $1$-strongly convex with respect to the product norm.

Run Algorithm~\ref{alg:compact-ftrl} with the norm \eqref{eq:compact-product-norm} and regularizer \eqref{eq:compact-product-potential}. The cone is closed, and the regularized objective is coercive and strictly convex, so each minimizer exists and is unique.

To bound regret, use Theorem~\ref{thm:ftrl-equality} with $\psi_{t}\coloneqq\psi$. Apply Lemma~\ref{lem:minimum-growth} with $p=F_{t}$ and $r=0$ in the product norm \eqref{eq:compact-product-norm}; the dual-norm calculation from Section~\ref{sec:bounded-euclidean} bounds each per-round cost.
Consequently, for every $\widecheck{\bu}\in\mathcal{C}_{\mathcal{V}}$,
\[
\sum_{t=1}^{T} \langle\widecheck{\bg}_{t},\widecheck{\by}_{t}-\widecheck{\bu}\rangle
\leqslant \psi(\widecheck{\bu}) + \frac{1}{2} \sum_{t=1}^{T} \|\widecheck{\bg}_{t}\|_{\star}^{2}\,.
\]
As in Section~\ref{sec:bounded-euclidean}, for $\bu\in\mathcal{V}$, define $\widecheck{\bu}\coloneqq(\bu,1)$ and compare with $\lambda\widecheck{\bu}$ for $\lambda>0$, which gives
\begin{equation}
\label{eq:compact-scale-optimization}
 \Regret_{T}(\bu)
 \leqslant\inf_{\lambda>0}\left\{\lambda\psi(\widecheck{\bu})+\frac{1}{2\lambda}\sum_{t=1}^{T}\|\widecheck{\bg}_{t}\|_{\star}^{2}\right\}
 =\sqrt{2\psi(\widecheck{\bu}) \sum_{t=1}^{T} \|\widecheck{\bg}_{t}\|_{\star}^{2}}\,.
\end{equation}
The comparator term satisfies
\[
2\psi(\widecheck{\bu})
=2\psi((\bu,1))
=\|\bu\|_{\mathcal{X}}^{2}+R^{2}
\leqslant2R^{2}\,.
\]
Also, by the dual norm definition,
\[
 \|\widecheck{\bg}_{t}\|_{\star}^{2}
 =\left\|(\bg_{t},-\langle \bg_{t},\bx_{t}\rangle)\right\|_{\star}^{2}
 =\frac{1}{\mu}\|\bg_{t}\|_{\mathcal{X},\star}^{2}+\frac{\langle \bg_{t},\bx_{t}\rangle^{2}}{R^{2}}
 \leqslant\left(1+\frac{1}{\mu}\right)\|\bg_{t}\|_{\mathcal{X},\star}^{2}\,.
\]
Combining these bounds, we obtain the following theorem.

\begin{theorem}
Let $\mathcal{X}$ be a finite-dimensional real normed space such that $\|\cdot\|_{\mathcal{X}}^{2}/2$ is $\mu$-strongly convex with respect to $\|\cdot\|_{\mathcal{X}}$, where $\mu>0$.
Let $\mathcal{V}\subseteq\mathcal{X}$ be nonempty, closed, bounded, and convex, and let $R>0$ satisfy $\|\bx\|_{\mathcal{X}}\leqslant R$ for every $\bx\in\mathcal{V}$.
Then Algorithm~\ref{alg:compact-ftrl} satisfies, for every $T\geqslant1$ and $\bu\in\mathcal{V}$,
\begin{equation}
\label{eq:compact-norm-regret}
 \Regret_{T}\left(\bu\right)
 \leqslant R\sqrt{2\left(1+\frac{1}{\mu}\right)}\sqrt{\sum_{t=1}^{T}\|\bg_{t}\|_{\mathcal{X},\star}^{2}}\,.
\end{equation}
\end{theorem}
This bound recovers the Euclidean bound when $\mu=1$.

For comparison, \ac{FTRL} directly on $\mathcal{V}$ with regularizer $\psi_{\mathcal{X}}/\eta$ gives the bound $\frac{R^{2}}{2\eta}+\frac{\eta}{2\mu}\sum_{t=1}^{T}\|\bg_{t}\|_{\mathcal{X},\star}^{2}$.
Its infimum over $\eta>0$ is $\frac{R}{\sqrt{\mu}}\sqrt{\sum_{t=1}^{T}\|\bg_{t}\|_{\mathcal{X},\star}^{2}}$, so our bound is a factor $\sqrt{2\left(1+\mu\right)}$ larger than this hindsight-optimized expression. However, the conic algorithm uses a fixed regularizer and requires no knowledge of the future subgradients.

\subsection{\texorpdfstring{\ac{CBA+}}{CBA+} with a strongly convex squared norm}
To recover the \ac{CBA+} of \citet{GrandClementK21} and extend it to the non-Euclidean case, we use the same conic lifting with \ac{OMD} as the regret minimizer.

\begin{algorithm}[t]
\caption{Non-Euclidean Conic Blackwell Algorithm Plus}
\label{alg:compact-md}
\algrenewcommand{\algorithmicrequire}{\textbf{Input:}}
\begin{algorithmic}[1]
\Require A feasible set $\mathcal{V}\subseteq\mathcal{X}$, a radius $R>0$, and $\bx_{0}\in\mathcal{V}$
\State $\mathcal{C}_{\mathcal{V}}\coloneqq\{\lambda(\bx,1):\bx\in\mathcal{V},\ \lambda\geqslant0\}$
\State $\displaystyle\psi((\by,\bar{y}))\coloneqq\frac{1}{2}\left(\|\by\|_{\mathcal{X}}^{2}+R^{2}\bar{y}^{2}\right)$
\State $\widecheck{\by}_{1}=(\by_{1},\bar{y}_{1})\coloneqq(\bzero,0)$ and $\widecheck{\bxi}_{1}\coloneqq\bzero$
\For{$t=1,2,\dots$}
\State $\bx_{t} = \decode(\widecheck{\by}_{t})$
 \State Play $\bx_{t}$ and observe $\bg_{t}\in\partial\ell_{t}(\bx_{t})$
 \State $\widecheck{\bg}_{t}\coloneqq(\bg_{t},-\langle\bg_{t},\bx_{t}\rangle)$
 \State $\displaystyle\widecheck{\by}_{t+1}=(\by_{t+1},\bar{y}_{t+1})\coloneqq\argmin_{\widecheck{\bx}\in\mathcal{C}_{\mathcal{V}}}\left\{\langle\widecheck{\bg}_{t},\widecheck{\bx}\rangle+B_{\psi}^{\widecheck{\bxi}_{t}}(\widecheck{\bx};\widecheck{\by}_{t})\right\}$
 \State Choose $\widecheck{\bxi}_{t+1}\in\partial\psi(\widecheck{\by}_{t+1})$ with
 $\widecheck{\bxi}_{t}-\widecheck{\bg}_{t}-\widecheck{\bxi}_{t+1}\in\mathcal{N}_{\mathcal{C}_{\mathcal{V}}}(\widecheck{\by}_{t+1})$
\EndFor
\end{algorithmic}
\end{algorithm}
\begin{proposition}
\label{prop:compact-md}
Run Algorithm~\ref{alg:compact-md} in a finite-dimensional real normed space or a real Hilbert space with its Hilbert norm. Assume $\tfrac{1}{2}\|\cdot\|_{\mathcal{X}}^{2}$ is $\mu$-strongly convex with respect to $\|\cdot\|_{\mathcal{X}}$ for $\mu>0$.
For every $T\geqslant1$ and every $\bu\in\mathcal{V}$, this algorithm satisfies
\[
 \Regret_{T}\left(\bu\right)
 \leqslant R\sqrt{2\left(1+\frac{1}{\mu}\right)}\sqrt{\sum_{t=1}^{T}\|\bg_{t}\|_{\mathcal{X},\star}^{2}}\,.
\]
\end{proposition}
\begin{proof}
In the Hilbert case, the update is a projection onto a closed convex cone; in the finite-dimensional case, coercivity gives attainment.
Strong convexity ensures uniqueness, and continuity of $\psi$ ensures the optimality condition provides the required subgradient.
We first bound regret on the cone using the one-step argument of \citet[Lemma~6.10]{Orabona19}.
For every $\widecheck{\bu}\in\mathcal{C}_{\mathcal{V}}$, constrained optimality with the subgradients specified in Algorithm~\ref{alg:compact-md} and the three-point Bregman identity give
\[
 \langle\widecheck{\bg}_{t},\widecheck{\by}_{t}-\widecheck{\bu}\rangle
 \leqslant B_{\psi}^{\widecheck{\bxi}_{t}}(\widecheck{\bu};\widecheck{\by}_{t}) - B_{\psi}^{\widecheck{\bxi}_{t+1}}(\widecheck{\bu};\widecheck{\by}_{t+1}) + \frac{1}{2}\|\widecheck{\bg}_{t}\|_{\star}^{2}\,.
\]
This calculation also covers nonsmooth $\psi$ and Hilbert spaces.
Summing yields
\[
 \sum_{t=1}^{T} \langle\widecheck{\bg}_{t},\widecheck{\by}_{t}-\widecheck{\bu}\rangle
 \leqslant B_{\psi}^{\widecheck{\bxi}_{1}}(\widecheck{\bu};\bzero) + \frac{1}{2}\sum_{t=1}^{T}\|\widecheck{\bg}_{t}\|_{\star}^{2}
 =\psi(\widecheck{\bu}) + \frac{1}{2}\sum_{t=1}^{T}\|\widecheck{\bg}_{t}\|_{\star}^{2}\,.
\]
For $\widecheck{\bu}\coloneqq(\bu,1)$, compare with $\lambda\widecheck{\bu}$ for $\lambda>0$.
Since $\langle \widecheck{\bg}_{t}, \widecheck{\by}_{t}\rangle=0$ and $\psi(\lambda\widecheck{\bu})=\lambda^{2}\psi(\widecheck{\bu})$, dividing by $\lambda$ and taking the infimum over $\lambda>0$ gives \eqref{eq:compact-scale-optimization}.
The comparator and lifted-gradient bounds used in \eqref{eq:compact-norm-regret} then give the stated guarantee.
\end{proof}

\subsection{Computing the update}
In Hilbert geometry, the update is projected gradient descent in the weighted norm $\|\cdot\|_{R}$:
\[
 \widecheck{\by}_{t+1}
 =\Pi_{\mathcal{C}_{\mathcal{V}}}^{R} \!\left(\widecheck{\by}_{t}- \left(\bg_{t},-\frac{\langle\bg_{t},\bx_{t}\rangle}{R^{2}}\right)\right)\,.
\]

\section{Closed convex constraints with a known gradient bound}
\label{sec:unbounded}

We now allow unbounded closed convex feasible sets, for which no finite radius $R$ is available. Their lifted cones need not be closed; Appendix~\ref{sec:open_cone} gives an example. We translate the domain by an anchor and replace the quadratic lifted regularizer by one with a logarithmic barrier. The barrier restricts its domain to positive height, where the lifted cone is relatively closed.

\citet{Kwon26approachability} shows how to use conic
lifting for unconstrained regret minimization: the closure of the lifted cone $\mathcal{C}_{\mathbb{R}^{d}}=\left(\mathbb{R}^{d}\times\mathbb{R}_{+}^{*}\right)\cup\{(\bzero,0)\}$ is $\mathbb{R}^{d}\times\mathbb{R}_{+}$, whose nonzero boundary points have zero scalar coordinate and cannot be decoded. To avoid these points, regularizers $\psi$ whose convex conjugates $\psi^{\star}$ satisfy
\[
 \nabla\psi^{\star}\left(\mathbb{R}^{d}\times\mathbb{R}\right)
 \subseteq\mathcal{C}_{\mathbb{R}^{d}}\,
\]
are proposed. Since the \ac{FTRL} iterate is $\nabla\psi^{\star}$ evaluated at the negative cumulative lifted loss, zero height can occur only at the origin. We extend this approach to general closed convex domains, including unbounded ones; our logarithmic barrier ensures positive height on every round.

Let $\mathcal{X}$ be a finite-dimensional real normed space or a real Hilbert space with its Hilbert norm.
Let $\mathcal{V}\subseteq\mathcal{X}$ be any nonempty closed convex set, possibly unbounded (for instance $\cV = \R_{+}^{d}$).
Assume a known bound $\|\bg_{t}\|_{\mathcal{X},\star}\leqslant G$, where $G>0$.
Assume that $\|\cdot\|_{\mathcal{X}}^{2}/2$ is $\mu$-strongly convex with respect to $\|\cdot\|_{\mathcal{X}}$ for some $\mu>0$.
Fix an anchor $\ba\in\mathcal{V}$ and set $\mathcal{K}\coloneqq\mathcal{V}-\ba$, so $\bzero\in\mathcal{K}$.
Fix an anchor-regret budget $\epsilon>0$, measured in loss units: the guarantee below keeps regret at $\ba$ below $\epsilon$.
Note that a zero budget would force a deterministic learner on $\mathbb{R}$ to remain at its anchor, as the following lemma shows.
\begin{lemma}
\label{lem:anchor-budget}
Fix $a\in\R$ and use scalar predictions $x_{t}$ and gradients $g_{t}$.
A deterministic online algorithm on $\R$ that satisfies $\Regret_{T}(a)\leqslant0$ for every gradient sequence with $|g_{t}|\leqslant1$ and every $T$ must play $x_{t}=a$ on every such sequence.
\end{lemma}
\begin{proof}
Suppose some admissible gradient history first makes the algorithm choose $x_{t}\ne a$ at round $t$.
Keep the preceding gradients unchanged and choose $g_{t}\coloneqq\operatorname{sign}(x_{t}-a)$.
This leaves $x_{t}$ unchanged because it is chosen before observing $g_{t}$.
All earlier predictions equal $a$, so
\[
 \Regret_{t}(a)=\sum_{i=1}^{t} g_{i}(x_{i}-a)=|x_{t}-a|>0\,.
\]
This contradicts the assumed guarantee. Thus, no admissible history can cause a first departure from the anchor.
\end{proof}

We now run \ac{FTRL} in the cone $\mathcal{C}_{\mathcal{K}}\coloneqq\{\lambda(\bz,1):\lambda\geqslant0,\ \bz\in\mathcal{K}\}$.
For $a_{t}>0$ and $\widecheck{\by}=(\by,\bar{y})$ with $\bar{y}>0$, define
\begin{equation}
\label{eq:barrier-potential}
 \psi_{t}(\widecheck{\by})
 = \psi_{t}((\by,\bar{y}))
 \coloneqq\frac{a_{t}}{2\epsilon}\|\by\|_{\mathcal{X}}^{2}+\epsilon\left(\bar{y}-\ln\bar{y}\right)\,.
\end{equation}
Set $\psi_{t}\coloneqq+\infty$ when $\bar{y}\leqslant0$.
For each played decision $\bx_{t}\in\mathcal{V}$, set $\bx_{t}^{\circ}\coloneqq\bx_{t}-\ba\in\mathcal{K}$.
Define $\widecheck{\bg}_{t}\coloneqq(\bg_{t},-\langle\bg_{t},\bx_{t}^{\circ}\rangle)$.
Since $\bx_{t}^{\circ}=\decode(\widecheck{\by}_{t})$, we again have $\langle\widecheck{\bg}_{t},\widecheck{\by}_{t}\rangle=0$.
The resulting algorithm is Algorithm~\ref{alg:barrier-ftrl}.
Theorem~\ref{thm:barrier-unified} establishes the existence of the updates and positive height in every round.

\begin{algorithm}[t]
\caption{Conic \ac{FTRL} on an unbounded closed convex set and known gradient bound}
\label{alg:barrier-ftrl}
\algrenewcommand{\algorithmicrequire}{\textbf{Input:}}
\begin{algorithmic}[1]
\Require A feasible set $\mathcal{V}\subseteq\mathcal{X}$, an anchor $\ba\in\mathcal{V}$, $G,\epsilon>0$, and $\beta\geqslant1$
\State $\mathcal{K}\coloneqq\mathcal{V}-\ba$ and $\mathcal{C}_{\mathcal{K}}\coloneqq\{\lambda(\bz,1):\lambda\geqslant0,\ \bz\in\mathcal{K}\}$
\State $\mathcal{D}\coloneqq\mathcal{C}_{\mathcal{K}}\cap(\mathcal{X}\times\mathbb{R}_{+}^{*})$ and $a_{t}\coloneqq\beta G^{2}t$ for every $t\geqslant1$
\State $\widecheck{\by}_{1}=(\by_{1},\bar{y}_{1})\coloneqq(\bzero,1)$
\For{$t=1,2,\dots$}
 \State $\bx_{t}^{\circ}\coloneqq \decode(\widecheck{\by}_{t})$ and play $\bx_{t}\coloneqq\ba+\bx_{t}^{\circ}$ \label{step: decoding knot x}
 \State Observe $\bg_{t}\in\partial\ell_{t}(\bx_{t})$ with $\|\bg_{t}\|_{\mathcal{X},\star}\leqslant G$
 \State $\widecheck{\bg}_{t}\coloneqq(\bg_{t},-\langle\bg_{t},\bx_{t}^{\circ}\rangle)$
 \State $\displaystyle\widecheck{\by}_{t+1}=(\by_{t+1},\bar{y}_{t+1})\coloneqq\argmin_{\widecheck{\by}=(\by,\bar{y})\in\mathcal{D}}\left\{\frac{a_{t+1}}{2\epsilon}\|\by\|_{\mathcal{X}}^{2}+\epsilon\left(\bar{y}-\ln\bar{y}\right)+\sum_{i=1}^{t}\langle\widecheck{\bg}_{i},\widecheck{\by}\rangle\right\}$
\EndFor
\end{algorithmic}
\end{algorithm}

For this schedule, $\beta\geqslant1$ ensures feasibility by Lemma~\ref{lem:barrier-feasible}; the general regret guarantee in Theorem~\ref{thm:barrier-unified} requires $\beta\geqslant2$.
In the unconstrained Hilbert case, $\beta=1$ recovers \ac{KT} coin-betting with its separate guarantee; see \eqref{eq:barrier-kt} and Remark~\ref{rem:kt-sharper}.

\subsection{Regret guarantee}
We now prove a regret guarantee for Algorithm~\ref{alg:barrier-ftrl} in real Hilbert spaces and in finite-dimensional real normed spaces with strongly convex squared norms.

\begin{theorem}[Conic \ac{FTRL} on closed convex sets]
\label{thm:barrier-unified}
Run Algorithm~\ref{alg:barrier-ftrl} with $\beta\geqslant2$ and $a_{t}\coloneqq\beta G^{2}t$ in a finite-dimensional real normed space or a real Hilbert space with its Hilbert norm.
Suppose $\psi(\by)\coloneqq\|\by\|_{\mathcal{X}}^{2}/2$ is $\mu$-strongly convex with respect to $\|\cdot\|_{\mathcal{X}}$ for $\mu>0$.
All updates are well defined, all decisions are feasible, and $\bar{y}_{t}>0$ on every round.
For every $T\geqslant1$ and $\bu\in\mathcal{V}$, simultaneously,
\begin{equation}
\label{eq:barrier-unified-regret}
 \Regret_{T}\left(\bu\right)\leqslant \epsilon+G\|\bu-\ba\|_{\mathcal{X}}
 \sqrt{\beta\left(T+1\right)\ln\left(1+
 \frac{\beta G^{2}\left(T+1\right)^{1+\left(2+1/\mu\right)/\beta}\|\bu-\ba\|_{\mathcal{X}}^{2}}{\epsilon^{2}}\right)}\,.
\end{equation}
In particular, $\Regret_{T}(\ba)<\epsilon$.
\end{theorem}

The analysis proceeds through existence and positive height (Lemmas~\ref{lem:conic-minimum} and~\ref{lem:barrier-feasible}), per-round control (Lemmas~\ref{lem:conic-energy} and~\ref{lem:barrier-stability}), and comparator-scale optimization (Lemma~\ref{lem:comparator-scale}). The Hilbert implementation follows in the next computational subsection; the unconstrained specialization is \ac{KT}, as shown in \eqref{eq:barrier-kt}.

\subsection{Existence and per-round control}
Define
\begin{equation}
\label{eq:barrier-state}
 \widecheck{\btheta}_{t}
 =(\btheta_{t}, \bar{\theta}_{t})
 \coloneqq\sum_{i=1}^{t} \widecheck{\bg}_{i}
 = \left(\sum_{i=1}^{t} \bg_{i}, -\sum_{i=1}^{t}\langle\bg_{i},\bx_{i}^{\circ}\rangle\right),
\end{equation}
and
\[
F_{t}(\widecheck{\by})
= F_{t}((\by,\bar{y}))
\coloneqq \frac{a_{t}}{2\epsilon}\|\by\|_{\mathcal{X}}^{2}+\epsilon\left(\bar{y}-\ln\bar{y}\right)+\sum_{i=1}^{t-1}\langle\widecheck{\bg}_{i},\widecheck{\by}\rangle~.
\]

We first give an existence lemma that applies to both unbounded-domain constructions.

\begin{lemma}[Existence and uniqueness of the lifted minimum]
\label{lem:conic-minimum}
Let $\mathcal{X}$ be a finite-dimensional real normed space whose squared norm is strongly convex, or a real Hilbert space with its Hilbert norm.
Let $\mathcal{D}_{t}\subseteq\mathcal{X}\times\mathbb{R}_{+}^{*}$ be nonempty, convex, and closed relative to this space.
For $\epsilon,a_{t}>0$, $\gamma\geqslant0$, $c_{t-1},\bar{\theta}_{t-1}\in\R$, and $\btheta_{t-1}\in\mathcal{X}^{\star}$, assume
\[
 \epsilon+c_{t-1}+\bar{\theta}_{t-1}>0\,.
\]
Then, the objective on $\mathcal{D}_{t}$, with $\widecheck{\by}=(\by,\bar{y})$,
\[
 F_{t}(\widecheck{\by})
 \coloneqq\frac{a_{t}}{2\epsilon}\|\by\|_{\mathcal{X}}^{2}
 +\frac{\gamma}{2\bar{y}}\|\by\|_{\mathcal{X}}^{2}
 +\left(\epsilon+c_{t-1}+\bar{\theta}_{t-1}\right)\bar{y}
 -\epsilon\ln\bar{y}+\langle\btheta_{t-1},\by\rangle,
\]
has a unique minimizer $\widecheck{\by}_{t}=(\by_{t},\bar{y}_{t})\in\mathcal{D}_{t}$, with $\bar{y}_{t}>0$.
\end{lemma}
\begin{proof}
The perspective term is nonnegative, and the dual norm inequality gives
\[
 F_{t}(\widecheck{\by})
 \geqslant\frac{a_{t}}{4\epsilon}\|\by\|_{\mathcal{X}}^{2}
 +\left(\epsilon+c_{t-1}+\bar{\theta}_{t-1}\right)\bar{y}
 -\epsilon\ln\bar{y}
 -\frac{\epsilon\|\btheta_{t-1}\|_{\mathcal{X},\star}^{2}}{a_{t}}\,.
\]
Since $\epsilon+c_{t-1}+\bar{\theta}_{t-1}>0$, the scalar term on the right tends to $+\infty$ as $\bar{y}\to0^{+}$ or $\bar{y}\to+\infty$. Thus, every sublevel set is bounded and stays away from $\bar{y}=0$.
Relative closedness of $\mathcal{D}_{t}$ and the logarithmic barrier make the extension of $F_{t}$ by $+\infty$ outside $\mathcal{D}_{t}$ proper, convex, and lower semicontinuous.
Its sublevel sets are compact in finite dimensions; in a Hilbert space, convexity and lower semicontinuity give weak lower semicontinuity, so these bounded sublevel sets are weakly compact.
Consequently, the minimum is attained.
Strong convexity of the squared norm and strict convexity of $-\ln\bar{y}$ give uniqueness.
\end{proof}

\subsection{Proof of the regret bound}

We start our proof with a preliminary lemma on coordinate bounds and a positivity induction.

\begin{lemma}[Coordinate bounds and positivity]
\label{lem:barrier-feasible}
Let $\mathcal{X}$ be a finite-dimensional real normed space whose squared norm is strongly convex, or a real Hilbert space with its Hilbert norm. Let $\mathcal{K}\subseteq\mathcal{X}$ be closed and convex with $\bzero\in\mathcal{K}$. For $\epsilon,a_{t}>0$, $\bar{\theta}_{t-1}>-\epsilon$, and $\btheta_{t-1}\in\mathcal{X}^{\star}$, the minimization step of Algorithm~\ref{alg:barrier-ftrl} has a unique minimizer $\widecheck{\by}_{t}=(\by_{t},\bar{y}_{t})$, satisfying
\begin{equation}
\label{eq:barrier-coordinate-bounds}
\begin{gathered}
 \frac{\epsilon}{\epsilon+\bar{\theta}_{t-1}}
 \leqslant \bar{y}_{t}\leqslant\frac{\epsilon}{\epsilon+\bar{\theta}_{t-1}}\left(1+\frac{\|\btheta_{t-1}\|_{\mathcal{X},\star}^{2}}{4a_{t}}\right),\\
 \|\by_{t}\|_{\mathcal{X}}
 \leqslant\frac{\epsilon\|\btheta_{t-1}\|_{\mathcal{X},\star}}{a_{t}},
 \ \text{and}\
 \|\by_{t}/\bar{y}_{t}\|_{\mathcal{X}}
 \leqslant\frac{(\epsilon+\bar{\theta}_{t-1})\|\btheta_{t-1}\|_{\mathcal{X},\star}}{a_{t}}\,.
\end{gathered}
\end{equation}
Suppose $G>0$, $\beta\geqslant1$, $a_{t}\geqslant\beta G^{2}t$, and $\|\bg_{t}\|_{\mathcal{X},\star}\leqslant G$ on every round.
Then, in Algorithm~\ref{alg:barrier-ftrl} all iterates are well defined and
\begin{equation}
\label{eq:barrier-coefficient-positive}
\begin{gathered}
 \|\bx_{t}^{\circ}\|_{\mathcal{X}}\leqslant
 \frac{G(t-1)}{a_{t}}(\epsilon+\bar{\theta}_{t-1})\leqslant\frac{t-1}{\beta Gt}(\epsilon+\bar{\theta}_{t-1}),
 \quad\text{and}\\
 \epsilon+\bar{\theta}_{t}
 \geqslant\left(1-\frac{G^{2}(t-1)}{a_{t}}\right)(\epsilon+\bar{\theta}_{t-1})
 \geqslant\left(1-\frac{t-1}{\beta t}\right)(\epsilon+\bar{\theta}_{t-1})>0\,.
\end{gathered}
\end{equation}
\end{lemma}
\begin{proof}
Since $\mathcal{K}$ is closed and the map $(\by,\bar{y})\mapsto\by/\bar{y}$ is continuous on $\mathcal{X}\times\mathbb{R}_{+}^{*}$, the feasible set $\mathcal{D}$ is closed relative to this space.
Apply Lemma~\ref{lem:conic-minimum} with $\mathcal{D}_{t}=\mathcal{D}$, $\gamma=0$, and $c_{t-1}=0$ to obtain the unique minimizer.

To bound $\|\by_{t}\|_{\mathcal{X}}$, note that $\bzero\in\mathcal{K}$ makes $(\lambda\by_{t},\bar{y}_{t})$ feasible for $0\leqslant\lambda\leqslant1$.
The objective along this segment is minimized at $\lambda=1$, so its left derivative there is nonpositive.
The definition of dual norm also gives a lower bound for this derivative:
\begin{align*}
 \frac{a_{t}}{\epsilon}\|\by_{t}\|_{\mathcal{X}}^{2}+\langle\btheta_{t-1},\by_{t}\rangle
 &\geqslant\frac{a_{t}}{\epsilon}\|\by_{t}\|_{\mathcal{X}}^{2}-\|\btheta_{t-1}\|_{\mathcal{X},\star}\|\by_{t}\|_{\mathcal{X}}\\
 &=\frac{a_{t}}{\epsilon}\left(\|\by_{t}\|_{\mathcal{X}}-\frac{\epsilon\|\btheta_{t-1}\|_{\mathcal{X},\star}}{2a_{t}}\right)^{2}
 -\frac{\epsilon\|\btheta_{t-1}\|_{\mathcal{X},\star}^{2}}{4a_{t}}
 \geqslant-\frac{\epsilon\|\btheta_{t-1}\|_{\mathcal{X},\star}^{2}}{4a_{t}}\,.
\end{align*}
Combining the nonpositive derivative with the definition of dual norm gives
\[
 a_{t}\|\by_{t}\|_{\mathcal{X}}^{2}/\epsilon
 \leqslant\|\btheta_{t-1}\|_{\mathcal{X},\star}\|\by_{t}\|_{\mathcal{X}}.
\]
This proves the bound on $\|\by_{t}\|_{\mathcal{X}}$.
Also, $\lambda\widecheck{\by}_{t}$ is feasible for every $\lambda>0$.
Along this curve the derivative at $\lambda=1$ is zero, giving
\[
 \frac{a_{t}}{\epsilon}\|\by_{t}\|_{\mathcal{X}}^{2}+\langle\btheta_{t-1},\by_{t}\rangle+(\epsilon+\bar{\theta}_{t-1})\bar{y}_{t}=\epsilon\,.
\]
The bounds on the first two terms give both bounds on $\bar{y}_{t}$.
Dividing the bound on $\|\by_{t}\|_{\mathcal{X}}$ by $\bar{y}_{t}\geqslant \epsilon/(\epsilon+\bar{\theta}_{t-1})$ gives the bound on $\|\by_{t}/\bar{y}_{t}\|_{\mathcal{X}}$.

For the sequential claim, initially $\bar{\theta}_{0}=0$ and $\widecheck{\by}_{1}=(\bzero,1)$.
Suppose $\epsilon+\bar{\theta}_{t-1}>0$. Lemma~\ref{lem:conic-minimum} gives the minimizer.
The coordinate estimate above, together with $\|\btheta_{t-1}\|_{\mathcal{X},\star}\leqslant G(t-1)$ and $a_{t}\geqslant\beta G^{2}t$, gives the bound on $\|\bx_{t}^{\circ}\|_{\mathcal{X}}$ in \eqref{eq:barrier-coefficient-positive}.
Using $\|\bg_{t}\|_{\mathcal{X},\star}\leqslant G$ in $\bar{\theta}_{t}=\bar{\theta}_{t-1}-\langle\bg_{t},\bx_{t}^{\circ}\rangle$ then gives the lower bound on $\epsilon+\bar{\theta}_{t}$ in the same display.
Finally, $G^{2}(t-1)/a_{t}\leqslant(t-1)/(\beta t)<1$, so $\epsilon+\bar{\theta}_{t}>0$ and the induction continues.
\end{proof}

The next lemma gives the key per-round and comparator bounds for the proof, covering both the known-gradient-bound construction and the adaptive construction of the next section.
\begin{lemma}
\label{lem:conic-energy}
Let $\mathcal{X}$ be a real normed space such that $\phi(\by)\coloneqq\|\by\|_{\mathcal{X}}^{2}/2$ is $\mu$-strongly convex with respect to $\|\cdot\|_{\mathcal{X}}$, where $\mu>0$.
Fix $T\geqslant1$, $\epsilon>0$, and $\gamma\geqslant0$.
Let $(a_{t})_{t=1}^{T+1}$ be positive and nondecreasing, let $(c_{t})_{t=0}^{T}$ be nondecreasing with $c_{0}=0$, and let
$\mathcal{D}_{t+1}\subseteq\mathcal{D}_{t}\subseteq\mathcal{X}\times\mathbb{R}_{+}^{*}$ be nonempty convex sets with $(\bzero,1)\in\mathcal{D}_{1}$.
For $\widecheck{\by}=(\by,\bar{y})\in\mathcal{D}_{t}$, define
\begin{equation}
\label{eq:conic-energy-objective}
\begin{aligned}
 \psi_{t}(\widecheck{\by})
 &\coloneqq\frac{a_{t}}{2\epsilon}\|\by\|_{\mathcal{X}}^{2}
 +\frac{\gamma}{2\bar{y}}\|\by\|_{\mathcal{X}}^{2}
 +\epsilon(\bar{y}-\ln\bar{y})+c_{t-1}\bar{y},\\
 F_{t}(\widecheck{\by})&\coloneqq\psi_{t}(\widecheck{\by})+\langle\widecheck{\btheta}_{t-1},\widecheck{\by}\rangle\,.
\end{aligned}
\end{equation}
Assume that $\widecheck{\by}_{t}\in\argmin_{\widecheck{\by}\in\mathcal{D}_{t}}\ F_{t}(\widecheck{\by})$ exists for each $t$.
For gradients $\bg_{t}\in\mathcal{X}^{\star}$, set $\bx_{t}\coloneqq\by_{t}/\bar{y}_{t}$, $\widecheck{\btheta}_{0}\coloneqq(\bzero,0)$, and
$\widecheck{\btheta}_{t}\coloneqq\widecheck{\btheta}_{t-1}+\widecheck{\bg}_{t}$, where $\widecheck{\bg}_{t}\coloneqq(\bg_{t},-\langle\bg_{t},\bx_{t}\rangle)$.
Suppose the decomposition $\bg_{t}=\bg_{t}^{\mathrm{c}}+\bg_{t}^{\mathrm{r}}$ satisfies
\begin{equation}
\label{eq:conic-energy-correction}
 \langle\bg_{t}^{\mathrm{c}},\by_{t}\rangle<\epsilon
 \quad\text{and}\quad
 \|\bg_{t}^{\mathrm{r}}\|_{\mathcal{X},\star}^{2}\leqslant2\mu\gamma(c_{t}-c_{t-1})\,.
\end{equation}
Then, for all $t$ we have
\begin{equation}
\label{eq:conic-energy-cost}
 F_{t}(\widecheck{\by}_{t})-F_{t+1}(\widecheck{\by}_{t+1})
 \leqslant-\frac{a_{t}(a_{t+1}-a_{t})}{2\epsilon a_{t+1}}\|\by_{t}\|_{\mathcal{X}}^{2}
 +\frac{\epsilon\|\bg_{t}^{\mathrm{c}}\|_{\mathcal{X},\star}^{2}}{2\mu a_{t+1}}
 -\frac{a_{t}}{a_{t+1}}\langle\bg_{t}^{\mathrm{c}},\by_{t}\rangle
 -\epsilon\ln\left(1-\frac{\langle\bg_{t}^{\mathrm{c}},\by_{t}\rangle}{\epsilon}\right)\,.
\end{equation}
Moreover, if $\kappa\geqslant0$ satisfies
\begin{equation}
\label{eq:conic-energy-sum}
 \sum_{t=1}^{T}\left(F_{t}(\widecheck{\by}_{t})-F_{t+1}(\widecheck{\by}_{t+1})\right)\leqslant\epsilon\kappa,
\end{equation}
then every $\bu\in\mathcal{X}$ and $\lambda>0$ with $\lambda\widecheck{\bu}\in\mathcal{D}_{T+1}$, where $\widecheck{\bu}\coloneqq(\bu,1)$, satisfy
\begin{equation}
\label{eq:conic-energy-comparator}
 \Regret_{T}(\bu)\leqslant\epsilon+c_{T}+\frac{\gamma}{2}\|\bu\|_{\mathcal{X}}^{2}
 +\frac{a_{T+1}\lambda\|\bu\|_{\mathcal{X}}^{2}}{2\epsilon}
 +\frac{\epsilon}{\lambda}(\kappa-1-\ln\lambda)\,.
\end{equation}
\end{lemma}
\begin{proof}
We first prove \eqref{eq:conic-energy-cost} by combining a lower bound on $F_{t}(\widecheck{\by}_{t+1})-F_{t}(\widecheck{\by}_{t})$ with the exact increment $F_{t+1}(\widecheck{\by}_{t+1})-F_{t}(\widecheck{\by}_{t+1})$.
The term $\gamma\|\by\|_{\mathcal{X}}^{2}/(2\bar{y})$ will absorb $\bg_{t}^{\mathrm{r}}$, and the logarithmic term will control the height ratio $\bar{y}_{t+1}/\bar{y}_{t}$.
We then use the assumed sum bound \eqref{eq:conic-energy-sum} to obtain the comparator inequality.

The nested domains give $\widecheck{\by}_{t+1}\in\mathcal{D}_{t+1}\subseteq\mathcal{D}_{t}$.
For $\widecheck{\by}=(\by,\bar{y})\in\mathcal{X}\times\mathbb{R}_{+}^{*}$, define the convex continuous function
\[
 r_{t}(\widecheck{\by})
 \coloneqq\frac{a_{t}}{2\epsilon}\|\by\|_{\mathcal{X}}^{2}+\frac{\gamma}{2\bar{y}}\|\by\|_{\mathcal{X}}^{2}-\epsilon\ln\bar{y}\,.
\]
On $\mathcal{D}_{t}$, the function $p\coloneqq F_{t}-r_{t}$ is affine.
Apply Lemma~\ref{lem:minimum-growth} with $f\coloneqq F_{t}$, $r\coloneqq r_{t}$, and the zero seminorm $|\widecheck{\by}|\coloneqq 0$.
The lemma supplies a subgradient $\widecheck{\bxi}_{t}\in\partial r_{t}(\widecheck{\by}_{t})$ for which \eqref{eq:minimum-growth} holds.
Holding $\bar{y}$ fixed shows that its first coordinate is $(a_{t}/\epsilon+\gamma/\bar{y}_{t})\bd_{t}$ for some $\bd_{t}\in\partial(\|\cdot\|_{\mathcal{X}}^{2}/2)(\by_{t})$; differentiating with respect to $\bar{y}$ gives its scalar coordinate. Thus,
\[
 \widecheck{\bxi}_{t}=(\bxi_{t},\bar{\xi}_{t})
 =\left(\left(\frac{a_{t}}{\epsilon}+\frac{\gamma}{\bar{y}_{t}}\right)\bd_{t},
 -\frac{\gamma\|\by_{t}\|_{\mathcal{X}}^{2}}{2\bar{y}_{t}^{2}}-\frac{\epsilon}{\bar{y}_{t}}\right)\,.
\]
Thus, \eqref{eq:minimum-growth} applied directly at $\widecheck{\by}_{t+1}$ gives
\[
\begin{aligned}
 F_{t}(\widecheck{\by}_{t+1})-F_{t}(\widecheck{\by}_{t})
 &\geqslant r_{t}(\widecheck{\by}_{t+1})-r_{t}(\widecheck{\by}_{t})-\langle\widecheck{\bxi}_{t},\widecheck{\by}_{t+1}-\widecheck{\by}_{t}\rangle\\
 &=\frac{a_{t}}{\epsilon}\left(\frac{\|\by_{t+1}\|_{\mathcal{X}}^{2}}{2}-\frac{\|\by_{t}\|_{\mathcal{X}}^{2}}{2}-\langle\bd_{t},\by_{t+1}-\by_{t}\rangle\right)\\
 &\quad+\frac{\gamma}{2\bar{y}_{t+1}}\|\by_{t+1}\|_{\mathcal{X}}^{2}
-\frac{\gamma}{2\bar{y}_{t}}\|\by_{t}\|_{\mathcal{X}}^{2}
-\frac{\gamma}{\bar{y}_{t}}
 \langle\bd_{t},\by_{t+1}-\by_{t}\rangle
+\frac{\gamma\|\by_{t}\|_{\mathcal{X}}^{2}}{2\bar{y}_{t}^{2}}
 (\bar{y}_{t+1}-\bar{y}_{t})\\
 &\quad+\epsilon\left(\frac{\bar{y}_{t+1}}{\bar{y}_{t}}-1-\ln\frac{\bar{y}_{t+1}}{\bar{y}_{t}}\right)\,.
\end{aligned}
\]
Homogeneity gives $(\bar{y}_{t+1}/\bar{y}_{t})\bd_{t}\in\partial(\|\cdot\|_{\mathcal{X}}^{2}/2)((\bar{y}_{t+1}/\bar{y}_{t})\by_{t})$, and Lemma~\ref{lem:squared-norm-subgradient} gives $\langle\bd_{t},\by_{t}\rangle=\|\by_{t}\|_{\mathcal{X}}^{2}$.
Consequently, strong convexity at $(\bar{y}_{t+1}/\bar{y}_{t})\by_{t}$ gives the following lower bound for the four terms proportional to $\gamma$:
\begin{align*}
\frac{\gamma}{2\bar{y}_{t+1}}&\|\by_{t+1}\|_{\mathcal{X}}^{2}
 -\frac{\gamma}{2\bar{y}_{t}}\|\by_{t}\|_{\mathcal{X}}^{2}
 -\frac{\gamma}{\bar{y}_{t}} \langle\bd_{t},\by_{t+1}-\by_{t}\rangle
 +\frac{\gamma\|\by_{t}\|_{\mathcal{X}}^{2}}{2\bar{y}_{t}^{2}} (\bar{y}_{t+1}-\bar{y}_{t})\\
&=
\frac{\gamma}{2\bar{y}_{t+1}}\|\by_{t+1}\|_{\mathcal{X}}^{2}
 -\frac{\gamma}{\bar{y}_{t}}\langle\bd_{t},\by_{t+1}\rangle
 +\frac{\gamma\bar{y}_{t+1}}{2\bar{y}_{t}^{2}}\|\by_{t}\|_{\mathcal{X}}^{2}\\
&=
\frac{\gamma}{\bar{y}_{t+1}}\left(\frac{\|\by_{t+1}\|_{\mathcal{X}}^{2}}{2}+\frac{(\bar{y}_{t+1})^{2}}{\bar{y}_{t}^{2}}\frac{\|\by_{t}\|_{\mathcal{X}}^{2}}{2}
 -\frac{\bar{y}_{t+1}}{\bar{y}_{t}}\langle\bd_{t},\by_{t+1}\rangle\right)\\
&=\frac{\gamma}{\bar{y}_{t+1}}\left(\frac{\|\by_{t+1}\|_{\mathcal{X}}^{2}}{2}-\frac{1}{2}\left\|\frac{\bar{y}_{t+1}}{\bar{y}_{t}}\by_{t}\right\|_{\mathcal{X}}^{2}
 -\frac{\bar{y}_{t+1}}{\bar{y}_{t}}\left\langle\bd_{t},\by_{t+1}-\frac{\bar{y}_{t+1}}{\bar{y}_{t}}\by_{t}\right\rangle\right)\\
 &\geqslant\frac{\mu\gamma}{2\bar{y}_{t+1}}\|\by_{t+1}-\bar{y}_{t+1}\bx_{t}\|_{\mathcal{X}}^{2}\,.
\end{align*}
Substituting this bound and using again $\langle\bd_{t},\by_{t}\rangle=\|\by_{t}\|_{\mathcal{X}}^{2}$ from Lemma~\ref{lem:squared-norm-subgradient} gives
\begin{equation}
\label{eq:energy-growth}
\begin{aligned}
 F_{t}(\widecheck{\by}_{t+1})-F_{t}(\widecheck{\by}_{t})
 &\geqslant\frac{a_{t}}{\epsilon}\left(\frac{\|\by_{t+1}\|_{\mathcal{X}}^{2}}{2}+\frac{\|\by_{t}\|_{\mathcal{X}}^{2}}{2}-\langle\bd_{t},\by_{t+1}\rangle\right)\\
 &\quad+\frac{\mu\gamma}{2\bar{y}_{t+1}}\|\by_{t+1}-\bar{y}_{t+1}\bx_{t}\|_{\mathcal{X}}^{2}\\
 &\quad+\epsilon\left(\frac{\bar{y}_{t+1}}{\bar{y}_{t}}-1-\ln\frac{\bar{y}_{t+1}}{\bar{y}_{t}}\right).
\end{aligned}
\end{equation}

The cumulative feedback and the coefficients give the exact objective increment
\[
 F_{t+1}(\widecheck{\by}_{t+1})-F_{t}(\widecheck{\by}_{t+1})
 =\frac{a_{t+1}-a_{t}}{2\epsilon}\|\by_{t+1}\|_{\mathcal{X}}^{2}
 +\langle\bg_{t},\by_{t+1}-\bar{y}_{t+1}\bx_{t}\rangle
 +(c_{t}-c_{t-1})\bar{y}_{t+1}\,.
\]
To absorb the contribution of $\bg_{t}^{\mathrm{r}}$ with the perspective term and $(c_{t}-c_{t-1})\bar{y}_{t+1}$, use Fenchel--Young and \eqref{eq:conic-energy-correction}. For $\gamma>0$, they give
\begin{equation}
\label{eq:surprise-completion}
 -\langle\bg_{t}^{\mathrm{r}},\by_{t+1}-\bar{y}_{t+1}\bx_{t}\rangle
 -\frac{\mu\gamma}{2\bar{y}_{t+1}}\|\by_{t+1}-\bar{y}_{t+1}\bx_{t}\|_{\mathcal{X}}^{2}
 \leqslant\frac{\bar{y}_{t+1}\|\bg_{t}^{\mathrm{r}}\|_{\mathcal{X},\star}^{2}}{2\mu\gamma}
 \leqslant(c_{t}-c_{t-1})\bar{y}_{t+1}\,.
\end{equation}
For $\gamma=0$, condition \eqref{eq:conic-energy-correction} forces $\bg_{t}^{\mathrm{r}}=\bzero$, and $c_{t}-c_{t-1}\geqslant0$ gives the same cancellation.
Combining this bound with \eqref{eq:energy-growth} and the objective increment leaves
\[
\begin{aligned}
 F_{t}(\widecheck{\by}_{t})-F_{t+1}(\widecheck{\by}_{t+1})
 &\leqslant\frac{a_{t}\langle\bd_{t},\by_{t+1}\rangle}{\epsilon}
 -\frac{a_{t+1}\|\by_{t+1}\|_{\mathcal{X}}^{2}+a_{t}\|\by_{t}\|_{\mathcal{X}}^{2}}{2\epsilon}\\
 &\quad-\langle\bg_{t}^{\mathrm{c}},\by_{t+1}-\bar{y}_{t+1}\bx_{t}\rangle
 -\epsilon\left(\frac{\bar{y}_{t+1}}{\bar{y}_{t}}-1-\ln\frac{\bar{y}_{t+1}}{\bar{y}_{t}}\right)\,.
\end{aligned}
\]
Split the remaining linear term as
\[
 -\langle\bg_{t}^{\mathrm{c}},\by_{t+1}-\bar{y}_{t+1}\bx_{t}\rangle
 =\langle\bg_{t}^{\mathrm{c}},\by_{t}-\by_{t+1}\rangle
 +\left(\frac{\bar{y}_{t+1}}{\bar{y}_{t}}-1\right)\langle\bg_{t}^{\mathrm{c}},\by_{t}\rangle\,.
\]
For every $x<1$, the maximum below is attained at $\alpha=1/(1-x)$:
\[
 \max_{\alpha>0}\left\{(\alpha-1)x-\alpha+1+\ln\alpha\right\}=-x-\ln(1-x)\,.
\]
By \eqref{eq:conic-energy-correction}, $x=\langle\bg_{t}^{\mathrm{c}},\by_{t}\rangle/\epsilon<1$. Applying the identity with this value of $x$ bounds the terms involving $\bar{y}_{t+1}/\bar{y}_{t}$ and gives
\begin{equation}
\label{eq:energy-quadratic-cost}
\begin{aligned}
 F_{t}(\widecheck{\by}_{t})-F_{t+1}(\widecheck{\by}_{t+1})
 &\leqslant\frac{a_{t}\langle\bd_{t},\by_{t+1}\rangle}{\epsilon}
 -\frac{a_{t+1}\|\by_{t+1}\|_{\mathcal{X}}^{2}+a_{t}\|\by_{t}\|_{\mathcal{X}}^{2}}{2\epsilon}\\
 &\quad+\langle\bg_{t}^{\mathrm{c}},\by_{t}-\by_{t+1}\rangle
 -\langle\bg_{t}^{\mathrm{c}},\by_{t}\rangle\\
 &\quad-\epsilon\ln\left(1-\frac{\langle\bg_{t}^{\mathrm{c}},\by_{t}\rangle}{\epsilon}\right).
\end{aligned}
\end{equation}

To obtain the term involving $a_{t+1}-a_{t}$ in \eqref{eq:conic-energy-cost}, apply strong convexity at $(a_{t}/a_{t+1})\by_{t}$, with subgradient $(a_{t}/a_{t+1})\bd_{t}$.
Using $\langle\bd_{t},\by_{t}\rangle=\|\by_{t}\|_{\mathcal{X}}^{2}$ gives the expansion
\begin{equation}
\label{eq:energy-rescaling}
\begin{aligned}
 \frac{a_{t+1}}{2}&\|\by_{t+1}\|_{\mathcal{X}}^{2}-a_{t}\langle\bd_{t},\by_{t+1}\rangle+\frac{a_{t}}{2}\|\by_{t}\|_{\mathcal{X}}^{2}\\
 &=a_{t+1}\left(\frac{\|\by_{t+1}\|_{\mathcal{X}}^{2}}{2}
 -\frac{1}{2}\left\|\frac{a_{t}}{a_{t+1}}\by_{t}\right\|_{\mathcal{X}}^{2}
 -\left\langle\frac{a_{t}}{a_{t+1}}\bd_{t},\by_{t+1}-\frac{a_{t}}{a_{t+1}}\by_{t}\right\rangle\right)
 +\frac{a_{t}(a_{t+1}-a_{t})}{2a_{t+1}}\|\by_{t}\|_{\mathcal{X}}^{2}\\
 &\geqslant\frac{a_{t}(a_{t+1}-a_{t})}{2a_{t+1}}\|\by_{t}\|_{\mathcal{X}}^{2}
 +\frac{\mu a_{t+1}}{2}\left\|\by_{t+1}-\frac{a_{t}\by_{t}}{a_{t+1}}\right\|_{\mathcal{X}}^{2}\,.
\end{aligned}
\end{equation}
The remaining linear term decomposes as
\[
 \langle\bg_{t}^{\mathrm{c}},\by_{t}-\by_{t+1}\rangle
 =\frac{a_{t+1}-a_{t}}{a_{t+1}}\langle\bg_{t}^{\mathrm{c}},\by_{t}\rangle
 -\left\langle\bg_{t}^{\mathrm{c}},\by_{t+1}-\frac{a_{t}}{a_{t+1}}\by_{t}\right\rangle\,.
\]
Combine this decomposition with \eqref{eq:energy-rescaling} divided by $\epsilon$. Fenchel--Young then gives
\[
 -\left\langle\bg_{t}^{\mathrm{c}},\by_{t+1}-\frac{a_{t}}{a_{t+1}}\by_{t}\right\rangle
 -\frac{\mu a_{t+1}}{2\epsilon}\left\|\by_{t+1}-\frac{a_{t}}{a_{t+1}}\by_{t}\right\|_{\mathcal{X}}^{2}
 \leqslant\frac{\epsilon\|\bg_{t}^{\mathrm{c}}\|_{\mathcal{X},\star}^{2}}{2\mu a_{t+1}}\,.
\]
Substitution into \eqref{eq:energy-quadratic-cost}, followed by combining the two coefficients of $\langle\bg_{t}^{\mathrm{c}},\by_{t}\rangle$, proves \eqref{eq:conic-energy-cost}.

For the regret bound, extend $\psi_{t}$ by $+\infty$ outside $\mathcal{D}_{t}$ and apply the \ac{FTRL} regret equality, Theorem~\ref{thm:ftrl-equality}, on $\mathcal{X}\times\mathbb{R}$ with comparator $\lambda\widecheck{\bu}$.
Since $c_{0}=0$, $\widecheck{\btheta}_{0}=(\bzero,0)$, and $(\bzero,1)\in\mathcal{D}_{1}$, the initial minimizer is $\widecheck{\by}_{1}=(\bzero,1)$ and $\psi_{1}(\widecheck{\by}_{1})=\epsilon$.
The initial condition is $\psi_{T+1}(\lambda\widecheck{\bu})-\epsilon$.
Since $\langle\widecheck{\bg}_{t},\widecheck{\by}_{t}\rangle=0$, the per-round terms in that equality are $F_{t}(\widecheck{\by}_{t})-F_{t+1}(\widecheck{\by}_{t+1})$; their sum is at most $\epsilon\kappa$ by \eqref{eq:conic-energy-sum}.
The final condition is nonpositive because $\lambda\widecheck{\bu}\in\mathcal{D}_{T+1}$.
The lifted regret equals $\lambda\Regret_{T}(\bu)$; expanding the terminal regularizer and dividing by $\lambda$ proves \eqref{eq:conic-energy-comparator}.
\end{proof}

The following lemma gives a sufficient condition for~\eqref{eq:conic-energy-sum}.
\begin{lemma}
\label{lem:barrier-stability}
Let $\mathcal{X}$ be a finite-dimensional real normed space or a real Hilbert space with its Hilbert norm,
and suppose $\psi(\by)\coloneqq\|\by\|_{\mathcal{X}}^{2}/2$ is $\mu$-strongly convex with respect to $\|\cdot\|_{\mathcal{X}}$ for $\mu>0$.
Let $\mathcal{K}\subseteq\mathcal{X}$ be closed and convex with $\bzero\in\mathcal{K}$,
and let $G,\epsilon>0$ and $\beta\geqslant2$.
Run Algorithm~\ref{alg:barrier-ftrl} on $\mathcal{V}\coloneqq\ba+\mathcal{K}$ with anchor $\ba$ and gradients satisfying $\|\bg_{t}\|_{\mathcal{X},\star}\leqslant G$.
Then, its \ac{FTRL} objective satisfies
\begin{equation}
\label{eq:potential-step}
 F_{t}(\widecheck{\by}_{t})-F_{t+1}(\widecheck{\by}_{t+1})
 \leqslant\frac{\epsilon}{\beta(t+1)}\left(1+\frac{1}{2\mu}\right)\,.
\end{equation}
\end{lemma}
\begin{proof}
Apply Lemma~\ref{lem:conic-energy} with $\mathcal{D}_{t}=\mathcal{D}$, $\gamma=0$, $c_{t}=0$, $\bg_{t}^{\mathrm{c}}=\bg_{t}$, $\bg_{t}^{\mathrm{r}}=\bzero$, and $a_{t}=\beta G^{2}t$, using the decoded decision $\bx_{t}^{\circ}$.
Lemma~\ref{lem:barrier-feasible} supplies the minimizers and gives $\|\by_{t}\|_{\mathcal{X}}\leqslant\epsilon(t-1)/(\beta Gt)$.
Consequently,
\[
 1-\max\left\{\frac{\langle\bg_{t},\by_{t}\rangle}{\epsilon},0\right\}
 \geqslant1-\frac{G\|\by_{t}\|_{\mathcal{X}}}{\epsilon}
 \geqslant\frac{(\beta-1)t+1}{\beta t}
 \geqslant\frac{t+1}{\beta t},
\]
where the last inequality uses $\beta\geqslant2$.
Applying \eqref{eq:logarithm-general-bound} with $x=\langle\bg_{t},\by_{t}\rangle/\epsilon$ therefore gives
\[
 -\epsilon\ln\left(1-\frac{\langle\bg_{t},\by_{t}\rangle}{\epsilon}\right)
 \leqslant\langle\bg_{t},\by_{t}\rangle
 +\frac{\langle\bg_{t},\by_{t}\rangle^{2}}{2\epsilon\left(1-\max\left\{\langle\bg_{t},\by_{t}\rangle/\epsilon,0\right\}\right)}
 \leqslant\langle\bg_{t},\by_{t}\rangle
 +\frac{\beta G^{2}t}{2\epsilon(t+1)}\|\by_{t}\|_{\mathcal{X}}^{2}\,.
\]
In \eqref{eq:conic-energy-cost}, this cancels the negative squared-norm term, since $a_{t}(a_{t+1}-a_{t})/a_{t+1}=\beta G^{2}t/(t+1)$.
Thus,
\[
 F_{t}(\widecheck{\by}_{t})-F_{t+1}(\widecheck{\by}_{t+1})
 \leqslant\frac{\langle\bg_{t},\by_{t}\rangle}{t+1}
 +\frac{\epsilon\|\bg_{t}\|_{\mathcal{X},\star}^{2}}{2\mu\beta G^{2}(t+1)}
 \leqslant\frac{\epsilon}{\beta(t+1)}\left(1+\frac{1}{2\mu}\right)\,.
\]
The last inequality uses $\|\bg_{t}\|_{\mathcal{X},\star}\leqslant G$ and $\|\by_{t}\|_{\mathcal{X}}\leqslant\epsilon/(\beta G)$.
\end{proof}

We now complete the proof by optimizing the comparator scale to obtain a parameter-free regret guarantee.
\begin{lemma}
\label{lem:comparator-scale}
Let $a,\epsilon>0$, $m,U,\kappa\geqslant0$, and $a\geqslant4m^{2}$, and define $\omega\coloneqq\ln(1+\frac{a e^{2\kappa}U^{2}}{\epsilon^{2}})$.
Then,
\begin{equation}
 \label{eq:comparator-scale}
 \inf_{\substack{\lambda>0\\2mU\lambda\leqslant\epsilon}}\left\{\frac{a\lambda U^{2}}{2\epsilon}+\frac{\epsilon}{\lambda}(\kappa-1-\ln\lambda)\right\}
 \leqslant U\max\left\{\sqrt{a\omega},2m\omega\right\}\,.
\end{equation}
\end{lemma}
\begin{proof}
If $U=0$, let $\lambda\to\infty$. Otherwise, $\omega>0$.
We bound the infimum by evaluating at an explicit feasible scale.
If $4m^{2}\omega\leqslant a$, take $\lambda\coloneqq\epsilon\sqrt{\omega}/(U\sqrt{a})$, which satisfies $2mU\lambda\leqslant\epsilon$.
Since $e^{\omega}-1=a e^{2\kappa}U^{2}/\epsilon^{2}$ and $1-e^{-\omega}\leqslant\omega$, the objective at this scale satisfies
\[
 \frac{a\lambda U^{2}}{2\epsilon}+\frac{\epsilon}{\lambda}\left(\kappa-1-\ln\lambda\right)
 =\frac{U\sqrt{a}}{2\sqrt{\omega}}\left(\omega+\ln\frac{e^{\omega}-1}{\omega}-2\right)
 \leqslant U\sqrt{a\omega}\,.
\]
Otherwise, $m>0$ and $a<4m^{2}\omega$.
Take $\lambda\coloneqq\epsilon/(2mU)$.
Using $a\geqslant4m^{2}$ and $\ln(a e^{2\kappa}U^{2}/\epsilon^{2})\leqslant\omega$, the value is at most
\[
 \frac{aU}{4m}+mU\left(2\kappa-2+\ln\frac{4m^{2}U^{2}}{\epsilon^{2}}\right)
 \leqslant mU\omega+mU(\omega-2)
 \leqslant2mU\omega\,. \qedhere
\]
\end{proof}

We are now ready to prove the regret guarantee.
\begin{proof}[Proof of Theorem~\ref{thm:barrier-unified}]
Lemma~\ref{lem:barrier-feasible} gives existence, uniqueness, $\bar{y}_{t}>0$, and $\epsilon+\bar{\theta}_{t}>0$.
Lemma~\ref{lem:barrier-stability} and $\sum_{t=1}^{T}1/(t+1)\leqslant\ln(T+1)$ give
\[
 \sum_{t=1}^{T}\left(F_{t}(\widecheck{\by}_{t})-F_{t+1}(\widecheck{\by}_{t+1})\right)
 \leqslant\frac{\epsilon}{\beta}\left(1+\frac{1}{2\mu}\right)\ln(T+1)\,.
\]
Apply \eqref{eq:conic-energy-comparator} with $\gamma=0$, $c_{t}=0$, $\mathcal{D}_{t}=\mathcal{D}$, $a_{t}=\beta G^{2}t$, $\kappa=(1+1/(2\mu))\ln(T+1)/\beta$, and centered comparator $\bu-\ba$.
Set $\widecheck{\bu}^{\circ}\coloneqq(\bu-\ba,1)$. Every $\lambda>0$ is feasible because $\lambda\widecheck{\bu}^{\circ}\in\mathcal{D}$.
Hence,
\begin{equation}
\label{eq:barrier-comparator-scale}
 \Regret_{T}(\bu)\leqslant \epsilon+\inf_{\lambda>0}\left\{
 \frac{\beta G^{2}(T+1)\lambda\|\bu-\ba\|_{\mathcal{X}}^{2}}{2\epsilon}
 +\frac{\epsilon}{\lambda}\left(\frac{1+1/(2\mu)}{\beta}\ln(T+1)-1-\ln\lambda\right)\right\}\,.
\end{equation}
Apply Lemma~\ref{lem:comparator-scale} with $a=\beta G^{2}(T+1)$, $m=0$, $U=\|\bu-\ba\|_{\mathcal{X}}$, and $\kappa=(1+1/(2\mu))\ln(T+1)/\beta$ to obtain \eqref{eq:barrier-unified-regret}.
Finally, $\epsilon+\bar{\theta}_{T}>0$ gives $\Regret_{T}(\ba)=-\bar{\theta}_{T}<\epsilon$.
\end{proof}

\subsection{Computing the update in the Hilbert case}
In a Hilbert space, the update reduces to scalar bisection, with one projection onto $\mathcal{K}$ per step and an explicit search interval. For the minimization step of Algorithm~\ref{alg:barrier-ftrl}, define
\begin{align}
 F_{t}(\widecheck{\by})
 &= F_{t}((\by,\bar{y}))
 \coloneqq\frac{a_{t}}{2\epsilon}\|\by\|_{2}^{2}+\epsilon\left(\bar{y}-\ln\bar{y}\right)+\langle\widecheck{\btheta}_{t-1},\widecheck{\by}\rangle \notag \\
 &= \frac{a_{t}}{2\epsilon}\left\|\by+\frac{\epsilon\btheta_{t-1}}{a_{t}}\right\|_{2}^{2} -\frac{\epsilon\|\btheta_{t-1}\|_{2}^{2}}{2a_{t}}+(\epsilon+\bar{\theta}_{t-1})\bar{y}-\epsilon\ln \bar{y}, \label{eq:barrier-objective}
\end{align}
where the last equality is by completing the square.
The minimizer can be computed by minimizing over $\by$ for each fixed $\bar{y}$, then minimizing over $\bar{y}$.
Write $\Pi_{\mathcal{K}}$ for the Hilbert-space projection onto $\mathcal{K}$, which returns the nearest point in $\mathcal{K}$.
Only the squared-distance term depends on $\by$, and the feasible values of $\by$ form $\bar{y}\mathcal{K}$.
Writing $\by=\bar{y}\bz$ with $\bz\in\mathcal{K}$ therefore gives
\begin{equation}
\label{eq:scalar-projection}
 \by_{t}(\bar{y})\coloneqq\argmin_{\by\in \bar{y}\mathcal{K}}\ \left\|\by+\frac{\epsilon\btheta_{t-1}}{a_{t}}\right\|_{2}^{2}
 =\bar{y}\,\Pi_{\mathcal{K}}\!\left(-\frac{\epsilon\btheta_{t-1}}{a_{t} \bar{y}}\right)\,.
\end{equation}
The remaining problem is one-dimensional:
\[
 \bar{y}_{t}=\argmin_{\bar{y}>0}\ F_{t}((\by_{t}(\bar{y}),\bar{y}))
 \quad\text{and}\quad \by_{t}=\by_{t}(\bar{y}_{t})\,.
\]
This objective is strictly convex: convex combinations of feasible minimizing pairs remain feasible, and $-\epsilon\ln \bar{y}$ makes the convexity inequality strict whenever the values of $\bar{y}$ differ.

To compute its derivative, use the fixed feasible set $\mathcal{K}$:
\[
 F_{t}((\by_{t}(\bar{y}),\bar{y}))=\min_{\bz\in\mathcal{K}}\left\{\frac{a_{t} \bar{y}^{2}}{2\epsilon}\|\bz\|_{2}^{2}
 +\bar{y}\langle\btheta_{t-1},\bz\rangle+(\epsilon+\bar{\theta}_{t-1})\bar{y}-\epsilon\ln \bar{y}\right\}\,.
\]
The unique minimizer is $\by_{t}(\bar{y})/\bar{y}$, which depends continuously on $\bar{y}>0$ because projection does not increase distances.
For two nearby values of $\bar{y}$, use the minimizer at each value as a trial point at the other to bound the change in the minimum from above and below.
After division by the difference, continuity makes both bounds converge to the derivative with $\bz$ fixed at $\by_{t}(\bar{y})/\bar{y}$.
Thus, no derivative of the projection is needed, and
\begin{equation}
\label{eq:scalar-projection-derivative}
 \frac{d}{d \bar{y}}F_{t}((\by_{t}(\bar{y}),\bar{y}))=\epsilon+\bar{\theta}_{t-1}-\frac{\epsilon}{\bar{y}}
 +\frac{1}{\bar{y}}\left(\frac{a_{t}}{\epsilon}\|\by_{t}(\bar{y})\|_{2}^{2}+\langle\btheta_{t-1},\by_{t}(\bar{y})\rangle\right)\,.
\end{equation}

The derivative is increasing and has a unique zero at $\bar{y}_{t}$.
By \eqref{eq:barrier-coordinate-bounds}, this zero lies in the explicit interval
\[
 \frac{\epsilon}{\epsilon+\bar{\theta}_{t-1}}\leqslant \bar{y}_{t}\leqslant
 \frac{\epsilon}{\epsilon+\bar{\theta}_{t-1}}\left(1+\frac{\|\btheta_{t-1}\|_{2}^{2}}{4a_{t}}\right)\,.
\]
Hence, bisection computes $\bar{y}_{t}$: let $\bar{y}$ be the arithmetic mean of the interval endpoints, then evaluate the projection in \eqref{eq:scalar-projection} and the derivative in \eqref{eq:scalar-projection-derivative}.
If the derivative is negative, replace the lower endpoint by $\bar{y}$; if it is positive, replace the upper endpoint; if it is zero, the minimizer has been found.
When $\btheta_{t-1}=\bzero$, the interval collapses to $\bar{y}_{t}=\epsilon/(\epsilon+\bar{\theta}_{t-1})$ and $\by_{t}=\bzero$.
Finally, play $\bx_{t}=\ba+\by_{t}(\bar{y}_{t})/\bar{y}_{t}$.
In finite dimensions, all other operations are vector arithmetic and inner products, so computing the projection onto $\mathcal{K}$ suffices to compute the update.

\subsection{The unconstrained update is \texorpdfstring{\ac{KT}}{KT} coin-betting}
When $\mathcal{V}=\mathcal{H}$ is a real Hilbert space, one can take $\beta=1$.
The minimizations over $\by$ and $\bar{y}$ in Algorithm~\ref{alg:barrier-ftrl} then separate, giving
\begin{equation}
\label{eq:barrier-kt}
 \by_{t}=-\frac{\epsilon\btheta_{t-1}}{G^{2}t},\qquad
 \bar{y}_{t}=\frac{\epsilon}{\epsilon+\bar{\theta}_{t-1}},\quad\text{and}\quad
 \bx_{t}^{\circ}=-\frac{(\epsilon+\bar{\theta}_{t-1})\btheta_{t-1}}{G^{2}t}\,.
\end{equation}
For comparison with \ac{KT}, define $\bg_{t}^{\mathrm{n}}\coloneqq\bg_{t}/G$.
Then, $(\epsilon+\bar{\theta}_{t-1})/G=\epsilon/G-\sum_{i=1}^{t-1}\langle\bg_{i}^{\mathrm{n}},\bx_{i}^{\circ}\rangle$
and $\btheta_{t-1}/G=\sum_{i=1}^{t-1}\bg_{i}^{\mathrm{n}}$.
Thus, \eqref{eq:barrier-kt} is exactly Algorithm~13.2 of \citet{Orabona19} in one dimension
and Algorithm~1 of \citet{OrabonaP16} in a Hilbert space with reward vectors $-\bg_{t}^{\mathrm{n}}$ and initial wealth $\epsilon/G$.

\begin{remark}[A sharper known regret bound for \ac{KT}]
\label{rem:kt-sharper}
In the unconstrained Hilbert case, we can use $\beta=1$, and a different proof using the \ac{KT} potential gives the sharper bound
\begin{equation}
\label{eq:kt-guarantees}
 \Regret_{T}\left(\bu\right)\leqslant \epsilon\left(1-\frac{1}{e\sqrt{T}}\right)
 +G\|\bu-\ba\|_{2}\sqrt{T\ln\left(1+\frac{e^{2}G^{2}T^{2}\|\bu-\ba\|_{2}^{2}}{\epsilon^{2}}\right)}\,.
\end{equation}
See \citet[Theorem~13.4 and Lemma~13.6]{Orabona19} and the Hilbert-space reduction of
\citet[Theorem~3]{OrabonaP16}.
\end{remark}

\subsection{Related algorithms}
\citet{CutkoskyO18} give two black-box reductions for non-Euclidean norms and constrained sets.
Their Algorithm~2 and Theorem~2 separate magnitude and direction: a scalar learner predicts $z_{t}$, a unit-ball learner predicts $\by_{t}$, and the played decision is $z_{t}\by_{t}$. The scalar feedback is $\langle\bg_{t},\by_{t}\rangle$, while the direction learner receives $\bg_{t}$. Their regret bound is
\[
 \Regret_{T}(\bu)\leqslant R_{T}^{\mathrm{scalar}}(\|\bu\|)
 +\|\bu\|R_{T}^{\mathrm{ball}}(\bu/\|\bu\|),
\]
with $\bu/\|\bu\|\coloneqq\bzero$ when $\bu=\bzero$. 
They then use a one-dimensional, parameter-free algorithm for the magnitude and a standard \ac{OCO} algorithm over balls for the direction.
Note that the scalar magnitude in their direction--magnitude reduction differs from our positive homogenizing coordinate.

Then, their Algorithm~3 and Theorem~3 transfer a learner $\mathcal{A}$ on a convex domain containing a nonempty closed convex set $\mathcal{W}$ in a reflexive Banach space to $\mathcal{W}$. Given its prediction $\bz_{t}$, play a metric projection onto $\mathcal{W}$ and return a subgradient at $\bz_{t}$ of
\[
 \widetilde{\ell}_{t}(\bx)
 \coloneqq\frac{\langle\bg_{t},\bx\rangle
 +\|\bg_{t}\|_{\star}\operatorname{dist}(\bx,\mathcal{W})}{2}\,.
\]
This gives $\Regret_{T}(\bu)\leqslant2R_{T}^{\mathcal{A}}(\bu)$ for every $\bu\in \mathcal{W}$, while preserving the gradient-norm bound. 
On the other hand, our approach directly constructs parameter-free algorithms for constrained sets with non-Euclidean norms.

\citet{Kwon26approachability} applies \ac{FTRL} to the
lifted loss vectors $\widecheck{\bg}_{t}$ over $\mathcal{C}_{\mathbb{R}^{d}}$,
with a fixed regularizer that, after normalizing the known gradient bound to one, takes the form
\[
 \psi((\by,\bar{y}))
 \coloneqq\frac{3}{2}\left(\bar{y}^{2}\varphi\left(\left\|\by/\bar{y}\right\|_{2}^{2}\right)+\bar{y}\chi\left(\left\|\by/\bar{y}\right\|_{2}\right)\right),
\]
for $\bar{y}>0$, where $\varphi(s)\coloneq1+s(1+\ln(1+s))$ and $\chi(s)\coloneqq(1+s)\ln(1+s)-s$ for $s\geqslant0$. The regularizer is zero at the origin and $+\infty$ elsewhere. Each prediction is a nonnegative multiple of the negative cumulative gradient, with its norm obtained by solving a scalar equation.

Let $\psi^{\star}$ denote the convex conjugate of $\psi$ and
$\widecheck{\btheta}_{T}$ denote this algorithm's cumulative lifted losses,
formed from its own predictions. The regret bound retains
$\psi^{\star}(-\widecheck{\btheta}_{T})$, which is at most
$\frac{1}{5}\sum_{t=1}^{T}\|\bg_{t}\|_{2}^{2}$ and can in some situations be much smaller. Optimizing the comparator scale gives a square-root dependence on this value and the additive comparator penalty $3\chi(\|\bu\|_{2})/2$. In particular, the stated regret bound at the origin is $\sqrt{\frac{6}{5}\sum_{t=1}^{T}\|\bg_{t}\|_{2}^{2}}$.

Our regularizer~\eqref{eq:barrier-potential} instead separates the
vector and scalar coordinates, increases its quadratic coefficient
with $t$, and enforces positive height through a logarithmic barrier. It
prescribes the anchor-regret budget $\epsilon$ and directly handles closed
convex constraints in Hilbert spaces and finite-dimensional spaces
with strongly convex squared norms. The guarantee in
Theorem~\ref{thm:barrier-unified} uses $G^{2}(T+1)$ rather than the
sum of observed squared gradient norms. In the unconstrained Hilbert
case with $\beta=1$, our update is explicitly \ac{KT} coin-betting, whereas the
fixed regularizer from \citet{Kwon26approachability} requires
solving a scalar equation to determine $\left\| \bx_{t} \right\|_{2}$.

\section{Unconstrained online optimization with an unknown gradient bound}
\label{sec:energy-adaptation}

We now set $\mathcal{V}=\mathcal{X}$ and, for simplicity, the anchor to $\bzero$, and remove the requirement that the algorithm receive a gradient bound.

Let $\mathcal{X}$ be a finite-dimensional real normed space, or a real Hilbert space with its Hilbert norm.
Assume that $\phi(\by)\coloneqq\|\by\|_{\mathcal{X}}^{2}/2$ is $\mu$-strongly convex with respect to $\|\cdot\|_{\mathcal{X}}$, where $0<\mu\leqslant1$.
For the Hilbert norm, take $\mu=1$.
The lifting itself is unchanged.

Fix $\epsilon,\gamma>0$ to specify an initial loss budget and a regularization tradeoff, without requiring a gradient bound, a comparator, or a horizon.

After observing $\bg_{t}\in\mathcal{X}^{\star}$, maintain
\begin{equation}
\label{eq:energy-statistics}
 m_{t}\coloneqq\max\left\{m_{t-1},\|\bg_{t}\|_{\mathcal{X},\star}\right\},\qquad
 v_{t}\coloneqq\sum_{i=1}^{t}\|\bg_{i}\|_{\mathcal{X},\star}^{2},
 \quad\text{and}\quad \delta_{t}\coloneqq\sum_{i=1}^{t}(m_{i}-m_{i-1})^{2}\,.
\end{equation}
Initialize $m_{0}\coloneq0$, $v_{0}\coloneq0$, and $\delta_{0}\coloneq0$.
Thus, $\delta_{t}$ sums the squared increases in the largest dual gradient norm observed so far.
The lifted cone is $\mathcal{C}_{\mathcal{X}}\coloneqq\left(\mathcal{X}\times\mathbb{R}_{+}^{*}\right)\cup\{(\bzero,0)\}$.
Before round $t$, set
\begin{equation}
\label{eq:energy-domain}
\begin{gathered}
 a_{t}\coloneqq\gamma\epsilon+4(v_{t-1}+m_{t-1}^{2}),\qquad
 \rho_{t}\coloneqq\begin{cases}
 \epsilon/(2m_{t-1}),&\text{if $m_{t-1}>0$},\\
 +\infty,&\text{otherwise},
 \end{cases}\\
 \text{and}\quad
 \mathcal{D}_{t}\coloneqq\{\widecheck{\by}=(\by,\bar{y})\in\mathcal{C}_{\mathcal{X}}:\bar{y}>0,\ \|\by\|_{\mathcal{X}}\leqslant\rho_{t}\}\,.
\end{gathered}
\end{equation}
The sets $\mathcal{D}_{t}$ are nested subsets of the homogenized cone. Once $\rho_{t}$ is finite, they are not cones: upward scaling can violate $\|\by\|_{\mathcal{X}}\leqslant\rho_{t}$. Nevertheless, every $\bx\in\mathcal{X}$ has a feasible representation $\bar{y}(\bx,1)$ for sufficiently small $\bar{y}>0$.
For $\widecheck{\by}\in\mathcal{D}_{t}$, define the regularizer
\begin{equation}
\label{eq:energy-regularizer}
 \psi_{t}(\widecheck{\by})
 \coloneqq\frac{a_{t}}{2\epsilon}\|\by\|_{\mathcal{X}}^{2} +\frac{\gamma}{2\bar{y}}\|\by\|_{\mathcal{X}}^{2}+\epsilon(\bar{y}-\ln\bar{y})+\frac{\delta_{t-1}}{2\mu\gamma}\bar{y}\,.
\end{equation}
The update at the end of round $t$ computes $\widecheck{\by}_{t+1}$ from $\widecheck{\btheta}_{t}$ and the statistics $m_{t},v_{t},\delta_{t}$.

Note that the four terms in \eqref{eq:energy-regularizer} have the following roles. The quadratic coefficient $a_{t}$ makes the per-round costs depend on the sum of squared dual gradient norms. The logarithmic barrier enforces positive height. The quadratic perspective controls feedback beyond the previous gradient scale, while the linear term involving $\delta_{t-1}$ compensates for increases in that scale. The constraint on $\|\by\|_{\mathcal{X}}$ bounds the inner products in the per-round estimates.

\begin{algorithm}[t]
\caption{Conic \ac{FTRL} with unknown gradient bound}
\label{alg:energy-ftrl}
\algrenewcommand{\algorithmicrequire}{\textbf{Input:}}
\begin{algorithmic}[1]
\Require A real normed space $\mathcal{X}$, finite-dimensional or Hilbert; $0<\mu\leqslant1$ such that $\|\cdot\|_{\mathcal{X}}^{2}/2$ is $\mu$-strongly convex in its norm; $\epsilon,\gamma>0$
\State $\mathcal{C}_{\mathcal{X}}\coloneqq\left(\mathcal{X}\times\mathbb{R}_{+}^{*}\right)\cup\left\{(\bzero,0)\right\}$
\State $\widecheck{\by}_{1}=(\by_{1},\bar{y}_{1})\coloneqq(\bzero,1)$ and $\widecheck{\btheta}_{0}\coloneqq(\bzero,0)$
\State $m_{0}\coloneq0$, $v_{0}\coloneq0$, and $\delta_{0}\coloneq0$
\For{$t=1,2,\dots$}
 \State Play $\bx_{t}\coloneqq\by_{t}/\bar{y}_{t}$ and observe $\bg_{t}\in\partial\ell_{t}(\bx_{t})$
 \State $\widecheck{\bg}_{t}\coloneqq(\bg_{t},-\langle\bg_{t},\bx_{t}\rangle)$ and $\widecheck{\btheta}_{t}\coloneqq\widecheck{\btheta}_{t-1}+\widecheck{\bg}_{t}$
 \State $m_{t}\coloneqq\max\left\{m_{t-1},\|\bg_{t}\|_{\mathcal{X},\star}\right\}$, $v_{t}\coloneqq v_{t-1}+\|\bg_{t}\|_{\mathcal{X},\star}^{2}$, and $\delta_{t}\coloneqq\delta_{t-1}+(m_{t}-m_{t-1})^{2}$
 \State $a_{t+1}\coloneqq\gamma\epsilon+4(v_{t}+m_{t}^{2})$
 \If{$m_{t}>0$}
  \State $\rho_{t+1}\coloneqq\epsilon/(2m_{t})$
 \Else
  \State $\rho_{t+1}\coloneqq+\infty$
 \EndIf
 \State $\mathcal{D}_{t+1}\coloneqq\left\{\widecheck{\by}=(\by,\bar{y})\in\mathcal{C}_{\mathcal{X}}:\bar{y}>0,\ \|\by\|_{\mathcal{X}}\leqslant\rho_{t+1}\right\}$
 \State $\displaystyle\psi_{t+1}(\widecheck{\by})\coloneqq\frac{a_{t+1}}{2\epsilon}\|\by\|_{\mathcal{X}}^{2}+\frac{\gamma}{2\bar{y}}\|\by\|_{\mathcal{X}}^{2}+\epsilon(\bar{y}-\ln\bar{y})+\frac{\delta_{t}}{2\mu\gamma}\bar{y}$ on $\mathcal{D}_{t+1}$
 \State $\displaystyle\widecheck{\by}_{t+1}=(\by_{t+1},\bar{y}_{t+1})\coloneqq\argmin_{\widecheck{\by}\in\mathcal{D}_{t+1}}\left\{\psi_{t+1}(\widecheck{\by})+\langle\widecheck{\btheta}_{t},\widecheck{\by}\rangle\right\}$
\EndFor
\end{algorithmic}
\end{algorithm}

\subsection{Regret guarantee}

\begin{theorem}
\label{thm:energy-ftrl}
Run Algorithm~\ref{alg:energy-ftrl}.
For every sequence of gradients, the \ac{FTRL} minimizer is unique and satisfies $\bar{y}_{t}>0$ on every round; initially, $\widecheck{\by}_{1}=(\bzero,1)$.
For every $T\geqslant1$ and $\bu\in\mathcal{X}$, define
\begin{equation}
\label{eq:energy-log}
 L_{T}(\bu)\coloneqq\ln\left(1+\frac{a_{T+1}\|\bu\|_{\mathcal{X}}^{2}}{\epsilon^{2}}
 \left(\frac{a_{T+1}}{\gamma\epsilon}\right)^{1+1/(4\mu)}\right)\,.
\end{equation}
Then, simultaneously for all such $T,\bu$,
\begin{equation}
\label{eq:energy-regret}
 \Regret_{T}(\bu)\leqslant\epsilon+\frac{\delta_{T}}{2\mu\gamma}
 +\frac{\gamma}{2}\|\bu\|_{\mathcal{X}}^{2}
 +\|\bu\|_{\mathcal{X}}\max\left\{\sqrt{a_{T+1}L_{T}(\bu)},\,2m_{T}L_{T}(\bu)\right\}\,.
\end{equation}
In particular,
\[
 \Regret_{T}(\bzero)
 \leqslant\epsilon+\frac{\delta_{T}}{2\mu\gamma}
 -\frac{\epsilon}{2\bar{y}_{T}}
 -\frac{(1-\mu)\gamma}{2}\|\bx_{T}\|_{\mathcal{X}}^{2}
 <\epsilon+\frac{\delta_{T}}{2\mu\gamma}\,.
\]
\end{theorem}

\begin{remark}
If $m_{t}$ increases from zero to $G$ in $k$ equal increments, with no further increases through round $T$, then $\delta_{T}=G^{2}/k$; a single jump gives $\delta_{T}=G^{2}$. Thus, the bound distinguishes gradual increases in gradient scale from a single jump.
\end{remark}

\subsection{Proof of the regret bound}

For $\bu\in\mathcal{X}$, define $\widecheck{\bu}\coloneqq(\bu,1)$.
As before, the conic lifting gives
\begin{equation}
\label{eq:energy-lifting}
 \langle\widecheck{\bg}_{t},\widecheck{\by}_{t}\rangle=0
 \quad\text{and}\quad
 \sum_{t=1}^{T}\langle\widecheck{\bg}_{t},\widecheck{\by}_{t}-\lambda\widecheck{\bu}\rangle
 =\lambda\Regret_{T}(\bu)\,.
\end{equation}
The constraint now restricts the comparator scale to $\lambda\|\bu\|_{\mathcal{X}}\leqslant\rho_{T+1}$. By Lemma~\ref{lem:comparator-scale}, this restriction gives the term $2m_{T}L_{T}(\bu)$ in the maximum in \eqref{eq:energy-regret}.

For the analysis, decompose the gradient as
\begin{equation}
\label{eq:energy-decomposition}
\begin{aligned}
 \bg_{t}^{\mathrm{c}}
 &\coloneqq\begin{cases}
 \bg_{t}\min\left\{1,m_{t-1}/\|\bg_{t}\|_{\mathcal{X},\star}\right\},&\text{if $\bg_{t}\ne\bzero$},\\
 \bzero,&\text{otherwise},
 \end{cases}\\
 \bg_{t}^{\mathrm{r}}
 &\coloneqq\bg_{t}-\bg_{t}^{\mathrm{c}}\,.
 \end{aligned}
\end{equation}
Then, $\|\bg_{t}^{\mathrm{c}}\|_{\mathcal{X},\star}\leqslant m_{t-1}$ and
$\|\bg_{t}^{\mathrm{r}}\|_{\mathcal{X},\star}^{2}=(m_{t}-m_{t-1})^{2}=\delta_{t}-\delta_{t-1}$.
This clipping is used only in the analysis; the algorithm uses the full lifted feedback. Clipping by running gradient norms and artificial expanding constraints already appear in \citet{Cutkosky19c}. \citet{MhammediK20} develop hinted scale-free algorithms and reductions that remove hints at an additional cost. Their lower bounds concern specified scale-homogeneous regret families and do not rule out the quadratic-in-gradient-scale tradeoff in \eqref{eq:energy-regret}.
Write the cumulative state and objective as
\begin{equation}
\label{eq:energy-state}
 \widecheck{\btheta}_{t}=(\btheta_{t},\bar{\theta}_{t})
 \coloneqq\left(\sum_{i=1}^{t}\bg_{i},-\sum_{i=1}^{t}\langle\bg_{i},\bx_{i}\rangle\right)
\end{equation}
and
\begin{equation}
\label{eq:energy-objective}
 F_{t}(\widecheck{\by})\coloneqq\psi_{t}(\widecheck{\by})
 +\langle\widecheck{\btheta}_{t-1},\widecheck{\by}\rangle\,.
\end{equation}

\begin{lemma}
\label{lem:energy-cost}
For every sequence of gradients, all minimizers in Algorithm~\ref{alg:energy-ftrl} exist uniquely and satisfy $\bar{y}_{t}>0$. With the state \eqref{eq:energy-state} and objective \eqref{eq:energy-objective}, we also have
\begin{equation}
\label{eq:energy-cost}
 \epsilon+\frac{\delta_{t}}{2\mu\gamma}+\bar{\theta}_{t}
 \geqslant\frac{\epsilon}{2\bar{y}_{t}}+\frac{(1-\mu)\gamma}{2}\|\bx_{t}\|_{\mathcal{X}}^{2}>0
 \quad\text{and}\quad
 F_{t}(\widecheck{\by}_{t})-F_{t+1}(\widecheck{\by}_{t+1})
 \leqslant\frac{(4\mu+1)\epsilon}{8\mu}\frac{a_{t+1}-a_{t}}{a_{t+1}}\,.
\end{equation}
\end{lemma}
\begin{proof}
We first prove existence, uniqueness, and the first inequality in \eqref{eq:energy-cost} by induction.
Initially, $\epsilon+\delta_{0}/(2\mu\gamma)+\bar{\theta}_{0}=\epsilon>0$ and $\widecheck{\by}_{1}=(\bzero,1)$.
Suppose $\epsilon+\delta_{t-1}/(2\mu\gamma)+\bar{\theta}_{t-1}>0$.
The domain $\mathcal{D}_{t}$ is nonempty, convex, and closed relative to $\mathcal{X}\times\mathbb{R}_{+}^{*}$, including when $\rho_{t}=+\infty$.
Apply Lemma~\ref{lem:conic-minimum} with $c_{t-1}=\delta_{t-1}/(2\mu\gamma)$ to obtain the unique minimizer.
Since  the constraint on $\by$ is independent of $\bar{y}$, the optimality condition for the scalar component gives
\begin{equation}
\label{eq:energy-scalar-stationarity}
 \epsilon+\frac{\delta_{t-1}}{2\mu\gamma}+\bar{\theta}_{t-1}
 =\frac{\epsilon}{\bar{y}_{t}}+\frac{\gamma}{2}\|\bx_{t}\|_{\mathcal{X}}^{2}\,.
\end{equation}
The constraint and clipping give $|\langle\bg_{t}^{\mathrm{c}},\by_{t}\rangle|\leqslant\epsilon/2$, including $m_{t-1}=0$.
Fenchel--Young's inequality with coefficient $\mu\gamma$ gives
\[
 \langle\bg_{t}^{\mathrm{r}},\bx_{t}\rangle
 \leqslant\frac{\mu\gamma}{2}\|\bx_{t}\|_{\mathcal{X}}^{2}
 +\frac{\|\bg_{t}^{\mathrm{r}}\|_{\mathcal{X},\star}^{2}}{2\mu\gamma}\,.
\]
Combining this bound on the residual-gradient terms with the state updates, \eqref{eq:energy-scalar-stationarity}, and the clipped inner-product bound yields
\begin{equation}
\label{eq:energy-coefficient-positive}
 \epsilon+\frac{\delta_{t}}{2\mu\gamma}+\bar{\theta}_{t}
 =\frac{\epsilon}{\bar{y}_{t}}\left(1-\frac{\langle\bg_{t}^{\mathrm{c}},\by_{t}\rangle}{\epsilon}\right)
 +\frac{\gamma}{2}\|\bx_{t}\|_{\mathcal{X}}^{2}-\langle\bg_{t}^{\mathrm{r}},\bx_{t}\rangle
 +\frac{\|\bg_{t}^{\mathrm{r}}\|_{\mathcal{X},\star}^{2}}{2\mu\gamma}
 \geqslant\frac{\epsilon}{2\bar{y}_{t}}+\frac{(1-\mu)\gamma}{2}\|\bx_{t}\|_{\mathcal{X}}^{2}>0\,.
\end{equation}
This closes the induction, proving positivity and the first inequality in \eqref{eq:energy-cost}, together with existence and uniqueness of all iterates.

For the per-round cost, apply Lemma~\ref{lem:conic-energy} with the domains of Algorithm~\ref{alg:energy-ftrl}, $c_{t}=\delta_{t}/(2\mu\gamma)$, and the decomposition \eqref{eq:energy-decomposition}.
The correction condition \eqref{eq:conic-energy-correction} holds with equality in its second inequality.
Since $|\langle\bg_{t}^{\mathrm{c}},\by_{t}\rangle|\leqslant\epsilon/2$, the logarithmic bound \eqref{eq:logarithm-upper-bound} turns \eqref{eq:conic-energy-cost} into
\begin{equation}
\label{eq:energy-cost-bound}
 F_{t}(\widecheck{\by}_{t})-F_{t+1}(\widecheck{\by}_{t+1})
 \leqslant-\frac{a_{t}(a_{t+1}-a_{t})}{2\epsilon a_{t+1}}\|\by_{t}\|_{\mathcal{X}}^{2}
 +\frac{\epsilon\|\bg_{t}^{\mathrm{c}}\|_{\mathcal{X},\star}^{2}}{2\mu a_{t+1}}
 +\frac{a_{t+1}-a_{t}}{a_{t+1}}\langle\bg_{t}^{\mathrm{c}},\by_{t}\rangle
 +\frac{\langle\bg_{t}^{\mathrm{c}},\by_{t}\rangle^{2}}{\epsilon}\,.
\end{equation}
The coefficient schedule and $\|\bg_{t}^{\mathrm{c}}\|_{\mathcal{X},\star}\leqslant m_{t-1}$ give
\[
\begin{gathered}
 a_{t}\geqslant4m_{t-1}^{2}\geqslant4\|\bg_{t}^{\mathrm{c}}\|_{\mathcal{X},\star}^{2}
 \quad\text{and}\\
 a_{t+1}-a_{t}=4\left(\|\bg_{t}\|_{\mathcal{X},\star}^{2}+m_{t}^{2}-m_{t-1}^{2}\right)
 \geqslant4\|\bg_{t}^{\mathrm{c}}\|_{\mathcal{X},\star}^{2}\,.
\end{gathered}
\]
Consequently,
\[
 \frac{a_{t}(a_{t+1}-a_{t})}{2a_{t+1}}
 \geqslant\frac{1}{4}\min\{a_{t},a_{t+1}-a_{t}\}
 \geqslant\|\bg_{t}^{\mathrm{c}}\|_{\mathcal{X},\star}^{2}\,.
\]
Thus, the dual norm inequality gives
\[
 \frac{\langle\bg_{t}^{\mathrm{c}},\by_{t}\rangle^{2}}{\epsilon}
 \leqslant\frac{\|\bg_{t}^{\mathrm{c}}\|_{\mathcal{X},\star}^{2}\|\by_{t}\|_{\mathcal{X}}^{2}}{\epsilon}
 \leqslant\frac{a_{t}(a_{t+1}-a_{t})}{2\epsilon a_{t+1}}\|\by_{t}\|_{\mathcal{X}}^{2},
\]
which cancels the negative squared-norm term in \eqref{eq:energy-cost-bound}.
For the two remaining terms, the coefficient increment and the clipping bound give
\[
\begin{gathered}
 \frac{\epsilon\|\bg_{t}^{\mathrm{c}}\|_{\mathcal{X},\star}^{2}}{2\mu a_{t+1}}
 \leqslant\frac{\epsilon}{8\mu}\frac{a_{t+1}-a_{t}}{a_{t+1}}
 \quad\text{and}\quad 
 \frac{a_{t+1}-a_{t}}{a_{t+1}}\langle\bg_{t}^{\mathrm{c}},\by_{t}\rangle
 \leqslant\frac{\epsilon}{2}\frac{a_{t+1}-a_{t}}{a_{t+1}}\,.
\end{gathered}
\]
Their sum is the claimed $(4\mu+1)\epsilon(a_{t+1}-a_{t})/(8\mu a_{t+1})$ bound.
\end{proof}

\begin{proof}[Proof of Theorem~\ref{thm:energy-ftrl}]
Lemma~\ref{lem:energy-cost} gives well-definedness and, since $\Regret_{T}(\bzero)=-\bar{\theta}_{T}$, the strict anchor bound.
To prove \eqref{eq:energy-regret}, we bound the sum of per-round costs and then optimize over the feasible comparator scales.
Summing the per-round estimate in Lemma~\ref{lem:energy-cost} and applying \eqref{eq:logarithm-ratio-bounds} with $x=a_{t+1}/a_{t}$ gives
\begin{equation}
\label{eq:energy-stability-sum}
 \sum_{t=1}^{T}\left(F_{t}(\widecheck{\by}_{t})-F_{t+1}(\widecheck{\by}_{t+1})\right)
 \leqslant\frac{(4\mu+1)\epsilon}{8\mu}\sum_{t=1}^{T}\frac{a_{t+1}-a_{t}}{a_{t+1}}
 \leqslant\frac{(4\mu+1)\epsilon}{8\mu}\ln\frac{a_{T+1}}{\gamma\epsilon}\,.
\end{equation}
Apply \eqref{eq:conic-energy-comparator} with $c_{t}=\delta_{t}/(2\mu\gamma)$, $\kappa=(4\mu+1)\ln(a_{T+1}/(\gamma\epsilon))/(8\mu)$, and comparator $\bu$.
The feasible scales satisfy $\lambda>0$ and $2m_{T}\lambda\|\bu\|_{\mathcal{X}}\leqslant\epsilon$.
This gives
\[
 \Regret_{T}(\bu)\leqslant\epsilon+\frac{\delta_{T}}{2\mu\gamma}+\frac{\gamma}{2}\|\bu\|_{\mathcal{X}}^{2}
 +\frac{a_{T+1}\lambda\|\bu\|_{\mathcal{X}}^{2}}{2\epsilon}
 +\frac{\epsilon}{\lambda}\left(\frac{4\mu+1}{8\mu}\ln\frac{a_{T+1}}{\gamma\epsilon}-1-\ln\lambda\right)\,.
\]
Lemma~\ref{lem:comparator-scale}, with $a=a_{T+1}$, $m=m_{T}$, $U=\|\bu\|_{\mathcal{X}}$, and $\kappa=(4\mu+1)\ln(a_{T+1}/(\gamma\epsilon))/(8\mu)$, gives \eqref{eq:energy-regret}, including $\bu=\bzero$ and $m_{T}=0$.
Its assumption $a\geqslant4m^{2}$ follows from $a_{T+1}=\gamma\epsilon+4(v_{T}+m_{T}^{2})$.
\end{proof}

\subsection{Computing the update}
In finite dimensions, the update reduces to maximization over the norm ball and a scalar root search.
If $\btheta_{t-1}\ne\bzero$, define $\bz_{t}$ and, for $p>0$, $\by_{t}(p)$ by
\[
 \bz_{t}\coloneqq\argmax_{\bz\in\mathcal{X}:\|\bz\|_{\mathcal{X}}\leqslant1}\ \langle-\btheta_{t-1},\bz\rangle
 \quad\text{and}\quad
 \by_{t}(p)\coloneqq\min\left\{\rho_{t},\frac{\epsilon\|\btheta_{t-1}\|_{\mathcal{X},\star}}{a_{t}+\gamma p}\right\}\bz_{t}\,.
\]
Strict convexity of the norm ball makes $\bz_{t}$ unique, with $\|\bz_{t}\|_{\mathcal{X}}=1$. By the dual norm inequality, $\by_{t}(p)$ minimizes the objective over $\by$ when $\bar{y}=\epsilon/p$.
The scalar stationarity condition \eqref{eq:energy-scalar-stationarity} determines $p$ from
\begin{equation}
\label{eq:energy-scalar}
 \epsilon+\frac{\delta_{t-1}}{2\mu\gamma}+\bar{\theta}_{t-1}
 =p+\frac{\gamma p^{2}}{2\epsilon^{2}}\|\by_{t}(p)\|_{\mathcal{X}}^{2}
 \quad\text{and}\quad \widecheck{\by}_{t}=(\by_{t}(p),\epsilon/p)\,.
\end{equation}
The right-hand side of the scalar equation in \eqref{eq:energy-scalar} is continuous and strictly increasing from zero to infinity: $p\|\by_{t}(p)\|_{\mathcal{X}}$ is nondecreasing, and the first term is $p$.
Hence, the root is unique and satisfies $0<p\leqslant\epsilon+\delta_{t-1}/(2\mu\gamma)+\bar{\theta}_{t-1}$.
If $\btheta_{t-1}=\bzero$, take $\by_{t}\coloneqq\bzero$ and $\bar{y}_{t}\coloneqq\epsilon/(\epsilon+\delta_{t-1}/(2\mu\gamma)+\bar{\theta}_{t-1})$.
In a Hilbert space, $\bz_{t}=-\btheta_{t-1}/\|\btheta_{t-1}\|_{2}$, so $\by_{t}(p)$ is the projection of $-\epsilon\btheta_{t-1}/(a_{t}+\gamma p)$ onto the ball of radius $\rho_{t}$.

A zero gradient leaves the entire state unchanged, so inserting zero-gradient rounds preserves subsequent predictions.
With the state before round $t$ fixed, the next prediction $\bx_{t+1}$ varies continuously with $\bg_{t}$.
For $\btheta_{t-1}\ne\bzero$, the unique maximizing direction is continuous by compactness in finite dimensions and by the explicit formula in a Hilbert space.
At $\btheta_{t-1}=\bzero$, the bound $\|\by_{t}(p)\|_{\mathcal{X}}\leqslant\epsilon\|\btheta_{t-1}\|_{\mathcal{X},\star}/a_{t}$ gives continuity of the vector component.
The positive scalar coefficient and the strictly increasing equation \eqref{eq:energy-scalar} give continuity of its root.
Together with the continuous state updates and the fact that a zero gradient leaves the state unchanged, these observations imply $\bx_{t+1}\to\bx_{t}$ as $\bg_{t}\to\bzero$, with the state before round $t$ fixed.

\subsection{Related algorithms}
\citet[Theorem~1]{CutkoskyM24} also obtain a regret bound without a supplied bound on future gradients, with a quadratic tradeoff between gradient scale and comparator norm. In \eqref{eq:energy-regret}, the corresponding terms are $\delta_{T}/(2\mu\gamma)$ and $\gamma\|\bu\|_{\mathcal{X}}^{2}/2$. Since $\delta_{T}\leqslant m_{T}^{2}$, gradual increases in record gradient norms can reduce the former. Our construction requires no initial gradient-scale estimate; theirs uses an initial estimate that need not bound future gradients. Our construction also applies to any norm whose square is strongly convex, whereas theirs uses the Euclidean norm.

For unconstrained Euclidean regret minimization without a supplied gradient bound, \citet{Kwon26approachability} applies \ac{FTRL} to the lifted loss vectors $\widecheck{\bg}_{t}$ over $\mathcal{C}_{\mathcal{X}}$, with $\mathcal{X}=\mathbb{R}^{d}$. With $\gamma>0$, the regularizer can be written as
\[
 \psi_{t}((\by,\bar{y}))
 \coloneqq\frac{4m_{t-1}^{2}}{\gamma}(1+\bar{y}\ln\bar{y})+\frac{\gamma}{2\bar{y}}\|\by\|_{2}^{2}+\bar{y}^{2}\varphi\left(\left\|\by/\bar{y}\right\|_{2}^{2}\right),
\]
for $\bar{y}>0$, where $\varphi(s)\coloneq1+s(1+\ln(1+s))$ for $s\geqslant0$. Set $\psi_{t}((\bzero,0))\coloneq4m_{t-1}^{2}/\gamma$ and extend it by $+\infty$ outside $\mathcal{C}_{\mathcal{X}}$. Only the coefficient $4m_{t-1}^{2}/\gamma$ of $1+\bar{y}\ln\bar{y}$ changes over time. Let $\psi_{t}^{\star}$ denote the convex conjugate of $\psi_{t}$, and let $\widecheck{\btheta}_{T}$ denote this algorithm's cumulative lifted losses, formed from its own predictions. The regret bound retains the conjugate value $\psi_{T+1}^{\star}(-\widecheck{\btheta}_{T})$, which is at most $v_{T}/2$ and in some situations is much smaller than $v_{T}$. It also includes the quadratic comparator penalty $\gamma\|\bu\|_{2}^{2}/2$ and terms depending on $m_{T}^{2}/\gamma$.

In Euclidean geometry, this regularizer and \eqref{eq:energy-regularizer} share the quadratic perspective $\gamma\|\by\|_{2}^{2}/(2\bar{y})$. Algorithm~\ref{alg:energy-ftrl} instead puts the cumulative squared gradient norms into $a_{t}$ and uses a logarithmic barrier to keep the scalar coordinate positive. The term $\delta_{T}/(2\mu\gamma)$ in \eqref{eq:energy-regret} compensates for record increases, distinguishing gradual increases from a single jump. The constraint $\|\by\|_{\mathcal{X}}\leqslant\rho_{t}$ restricts comparator scales and gives the term $2m_{T}L_{T}(\bu)$. The maximum norm still enters $a_{T+1}$ and $L_{T}(\bu)$. The construction also applies beyond Euclidean spaces under the strong convexity assumption on the squared norm stated at the start of this section.

\section{Conclusion}

We have presented a unifying framework to design learning-rate-free online convex optimization algorithms. The core idea is to use conic lifting and augmented subgradients to rescale the comparator norm for free in the regret analysis, obtaining the same effect as oracle tuning of the learning rate. Extending the approach of \citet{Kwon26approachability} to general closed convex domains and non-Euclidean geometries allows us to generalize and recover both Conic Blackwell algorithms and parameter-free ones.

This note has only scratched the surface of what conic lifting can offer to the optimization world. We believe that many other applications of this framework are possible. For example, a natural direction is to derive an \ac{OMD} version of Algorithm~\ref{alg:barrier-ftrl} and establish comparable parameter-free regret guarantees. 

\section*{Acknowledgments}
The authors acknowledge the use of GPT-6-Astra and GPT-6.1-Sol for assistance with the proof of Section~\ref{sec:energy-adaptation} and for help in writing. The authors have checked all details for correctness and take full responsibility for the content of this paper.
This work benefited from insightful discussions with Vianney Perchet and Sylvain Sorin. JGC was supported by Hi!\ Paris and a French Agence
Nationale de la Recherche (ANR) grant, reference
\texttt{11-LABX-0047}. JK was supported by a French Agence
Nationale de la Recherche (ANR) grant, reference
\texttt{ANR-21-CE40-0020} (CONVERGENCE project), and by a public grant
as part of the \emph{Investissement d'avenir project}, reference
\texttt{ANR-11-LABX-0056-LMH}, LabEx LMH.

\bibliographystyle{abbrvnat_nourl}
\bibliography{learning}

@article{d2021optimistic,
  title={Optimistic and adaptive lagrangian hedging},
  author={D'Orazio, Ryan and Huang, Ruitong},
  journal={arXiv preprint arXiv:2101.09603},
  year={2021}
}

@article{gordon2006no,
  title={No-regret algorithms for online convex programs},
  author={Gordon, Geoffrey J},
  journal={Advances in Neural Information Processing Systems},
  volume={19},
  year={2006}
}

@inproceedings{abernethy2011blackwell,
  title={{Blackwell} approachability and no-regret learning are equivalent},
  author={Abernethy, Jacob and Bartlett, Peter L and Hazan, Elad},
  booktitle={Proceedings of the 24th Annual Conference on Learning Theory},
  pages={27--46},
  year={2011},
  organization={JMLR Workshop and Conference Proceedings}
}

@article{farina2023regret,
  title={Regret {M}atching+:({I}n){S}tability and fast convergence in games},
  author={Farina, Gabriele and Grand-Cl{\'e}ment, Julien and Kroer, Christian and Lee, Chung-Wei and Luo, Haipeng},
  journal={Advances in Neural Information Processing Systems},
  volume={36},
  pages={61546--61572},
  year={2023}
}

@inproceedings{Cutkosky19c,
  author    = {Cutkosky, A.},
  editor    = {Beygelzimer, A. and Hsu, D.},
  title     = {Artificial Constraints and Hints for Unbounded Online Learning},
  booktitle = {Conference on Learning Theory, {COLT} 2019, 25-28 June 2019, Phoenix,
               AZ, {USA}},
  series    = {Proceedings of Machine Learning Research},
  volume    = {99},
  pages     = {874--894},
  publisher = {{PMLR}},
  year      = {2019},
  url       = {http://proceedings.mlr.press/v99/cutkosky19a.html},
}

@article{Cesa-BianchiCG04,
	author = {N. Cesa-Bianchi and A. Conconi and Gentile, C. },
	journal = {IEEE Trans. Inf. Theory},
	number = {9},
	pages = {2050--2057},
	title = {On the generalization ability of on-line learning algorithms},
	volume = {50},
	year = {2004},
	url = {https://homes.di.unimi.it/~cesabian/Pubblicazioni/J20.pdf}
}

@article{zhang2025scale,
  title={Scale-invariant {R}egret {M}atching and online learning with optimal convergence: Bridging theory and practice in zero-sum games},
  author={Zhang, Brian Hu and Anagnostides, Ioannis and Sandholm, Tuomas},
  journal={arXiv preprint arXiv:2510.04407},
  year={2025}
}

@article{soleymani2026multiplicative,
  title={Multiplicative Optimism for Constant Regret in Games},
  author={Soleymani, Ashkan and Piliouras, Georgios},
  journal={arXiv preprint arXiv:2609.21976},
  year={2026}
}

@book{nesterov1994interior,
  title={Interior-point polynomial algorithms in convex programming},
  author={Nesterov, Yurii and Nemirovskii, Arkadii},
  year={1994},
  publisher={SIAM}
}

@article{soleymani2025cautious,
  title={Cautious optimism: A meta-algorithm for near-constant regret in general games},
  author={Soleymani, Ashkan and Piliouras, Georgios and Farina, Gabriele},
  journal={arXiv preprint arXiv:2506.05005},
  year={2025}
}

@article{liu2026constant,
  title={Constant individual regret in general games},
  author={Liu, Mingyang and Farina, Gabriele and Ozdaglar, Asuman},
  journal={arXiv preprint arXiv:2608.31166},
  year={2026}
}

@article{abbadi2026constant,
  title={Constant regret in general games via higher-order optimism},
  author={Abbadi, Omar and Laraki, Rida and Mertikopoulos, Panayotis},
  journal={arXiv preprint arXiv:2609.04113},
  year={2026}
}

@article{farina2022near,
  title={Near-optimal no-regret learning dynamics for general convex games},
  author={Farina, Gabriele and Anagnostides, Ioannis and Luo, Haipeng and Lee, Chung-Wei and Kroer, Christian and Sandholm, Tuomas},
  journal={Advances in Neural Information Processing Systems},
  volume={35},
  pages={39076--39089},
  year={2022}
}

@inproceedings{luo2022adaptive,
  title={Adaptive bandit convex optimization with heterogeneous curvature},
  author={Luo, Haipeng and Zhang, Mengxiao and Zhao, Peng},
  booktitle={Conference on Learning Theory},
  pages={1576--1612},
  year={2022},
  organization={PMLR}
}

@article{lee2020bias,
  title={Bias no more: high-probability data-dependent regret bounds for adversarial bandits and {MDP}s},
  author={Lee, Chung-Wei and Luo, Haipeng and Wei, Chen-Yu and Zhang, Mengxiao},
  journal={Advances in neural information processing systems},
  volume={33},
  pages={15522--15533},
  year={2020}
}

@book{Cesa-BianchiL06,
  author    = {Cesa-Bianchi, N.  and Lugosi, G. },
  title     = {Prediction, learning, and games},
  publisher = {Cambridge University Press},
  year      = {2006},
  url = {https://www.cambridge.org/us/academic/subjects/computer-science/pattern-recognition-and-machine-learning/prediction-learning-and-games?format=HB&isbn=9780521841085}
}

@article{AuerCG02,
  author    = {P. Auer and N. Cesa-Bianchi and C. Gentile},
  title     = {Adaptive and Self-Confident On-Line Learning Algorithms},
  journal   = {J. Comput. Syst. Sci.},
  volume    = {64},
  number    = {1},
  year      = {2002},
  pages     = {48--75},
  url       = {http://homes.dsi.unimi.it/~cesabian/Pubblicazioni/jcss-02.pdf}
}

@incollection{StreeterM12,
title ={No-Regret Algorithms for Unconstrained Online Convex Optimization},
author={M. Streeter and B. McMahan},
booktitle = {Advances in Neural Information Processing Systems 25},
pages = {2402--2410},
year = {2012},
publisher = {Curran Associates, Inc.},
url = {https://arxiv.org/pdf/1211.2260.pdf}
}

@article{Shalev-Shwartz12,
  author    = {S. Shalev-Shwartz},
  title     = {Online Learning and Online Convex Optimization},
  journal   = {Foundations and Trends in Machine Learning},
  volume    = {4},
  number = {2},
  year      = {2012}
}

@inproceedings{Zinkevich03,
  author    = {Zinkevich, M.},
  title     = {Online Convex Programming and Generalized Infinitesimal Gradient Ascent},
  booktitle = {Proc. of the International Conference on Machine Learning},
  year      = {2003},
  pages     = {928--936},
  url = {http://ra.adm.cs.cmu.edu/anon/usr0/ftp/usr/anon/2003/CMU-CS-03-110.pdf}
}

@article{StreeterM10,
Author = {M. Streeter and H. B. McMahan},
Title = {Less Regret via Online Conditioning},
Year = {2010},
journal={arXiv preprint arXiv:1002.4862},
url = {https://arxiv.org/abs/1002.4862}
}

@inproceedings{Orabona14,
author = {F. Orabona},
title = {Simultaneous Model Selection and Optimization through Parameter-free Stochastic Learning},
booktitle = {Advances in Neural Information Processing Systems 27},
year = 2014,
url={https://arxiv.org/abs/1406.3816}
}

@inproceedings{ChaudhuriYH09,
  title={A parameter-free hedging algorithm},
  author={Chaudhuri, K. and Freund, Y. and Hsu, D. J.},
  booktitle={Advances in neural information processing systems},
  pages={297--305},
  year={2009},
  url = {https://arxiv.org/abs/0903.2851}
}

@incollection{OrabonaP16,
title = {Coin Betting and Parameter-Free Online Learning},
author = {Orabona, F. and P\'al, D.},
booktitle = {Advances in Neural Information Processing Systems 29},
editor = {D. D. Lee and M. Sugiyama and U. V. Luxburg and I. Guyon and R. Garnett},
pages = {577--585},
year = {2016},
publisher = {Curran Associates, Inc.},
url = {https://arxiv.org/pdf/1602.04128.pdf}
}

@incollection{vanErvenK16,
title = {{MetaGrad}: Multiple Learning Rates in Online Learning},
author = {van Erven, T. and Koolen, W. M.},
booktitle = {Advances in Neural Information Processing Systems 29},
editor = {D. D. Lee and M. Sugiyama and U. V. Luxburg and I. Guyon and R. Garnett},
pages = {3666--3674},
year = {2016},
publisher = {Curran Associates, Inc.},
url={https://proceedings.neurips.cc/paper_files/paper/2016/file/14cfdb59b5bda1fc245aadae15b1984a-Paper.pdf}
}

@inproceedings{CutkoskyO18,
        author = {Cutkosky, A. and Orabona, F.},
        title = {Black-Box Reductions for Parameter-free Online Learning in {Banach} Spaces},
        booktitle = {Proc. of the Conference on Learning Theory (COLT)},
        year = 2018,
        url = {https://arxiv.org/abs/1802.06293}
}

@inproceedings{KoolenVE15,
  title={Second-order quantile methods for experts and combinatorial games},
  author={Koolen, W. M. and van Erven, T.},
  booktitle={Proc. of the Conference On Learning Theory},
  pages={1155--1175},
  year={2015},
  url = {https://arxiv.org/abs/1502.08009}
}

@inproceedings{McMahanS10,
  title={Adaptive Bound Optimization for Online Convex Optimization},
  author={H. B. McMahan and M. J. Streeter},
  booktitle={COLT},
  year={2010},
  url = {https://static.googleusercontent.com/media/research.google.com/en//pubs/archive/36483.pdf}
}

@inproceedings{DuchiHS10,
  title = {Adaptive Subgradient Methods for Online Learning and Stochastic Optimization},
  author = {J. Duchi and E. Hazan and Y. Singer},
  booktitle={COLT},
  year={2010},
  url = {https://stanford.edu/~jduchi/projects/DuchiHaSi10_colt.pdf}
}

@inproceedings{OrabonaT17,
  title={Training deep networks without learning rates through coin betting},
  author={Orabona, F. and Tommasi, T.},
  booktitle={Advances in Neural Information Processing Systems},
  pages={2160--2170},
  year={2017},
  url = {https://arxiv.org/abs/1705.07795}
}

@inproceedings{LuoS15,
  title={Achieving all with no parameters: {AdaNormalHedge}},
  author={Luo, H. and Schapire, R. E.},
  booktitle={Proc. of the Conference on Learning Theory (COLT)},
  pages={1286--1304},
  year={2015},
  url = {http://proceedings.mlr.press/v40/Luo15.pdf}
}

@inproceedings{ChernovV10,
  title={Prediction with advice of unknown number of experts},
  author={Chernov, A. and Vovk, V.},
  booktitle={Proc. of the Conference on Uncertainty in Artificial Intelligence (UAI)},
  year={2010},
  url = {https://arxiv.org/abs/1408.2040}
}

@InProceedings{KempkaKW19,
  title = 	 {Adaptive Scale-Invariant Online Algorithms for Learning Linear Models},
  author = 	 {Kempka, M. and Kot{\l}owski, W. and Warmuth, M. K.},
  booktitle = 	 {Proc. of the 36th International Conference on Machine Learning},
  pages = 	 {3321--3330},
  year = 	 {2019},
  editor = 	 {Chaudhuri, K. and Salakhutdinov, R.},
  volume = 	 {97},
  series = 	 {Proc. of Machine Learning Research},
  address = 	 {Long Beach, California, USA},
  month = 	 {09--15 Jun},
  publisher = 	 {PMLR},
}

@article{Orabona19,
  title={Online Learning: A Modern Introduction Using Convex Optimization},
  author={Orabona, F.},
  journal={arXiv preprint arXiv:1912.13213},
  year={2026},
  note={Version 10},
  url={https://arxiv.org/abs/1912.13213}
}

@article{Hazan16,
url = {http://dx.doi.org/10.1561/2400000013},
year = {2016},
volume = {2},
journal = {Foundations and Trends in Optimization},
title = {Introduction to Online Convex Optimization},
number = {3--4},
pages = {157--325},
author = {Hazan, E.}
}

@inproceedings{Gordon99b,
  title={Regret bounds for prediction problems},
  author={Gordon, G. J.},
  booktitle={Proc. of the twelfth annual conference on Computational learning theory (COLT)},
  pages={29--40},
  year={1999},
  url = {http://www.cs.cmu.edu/~ggordon/colt99.ps.gz}
}

@inproceedings{MhammediK20,
  title={Lipschitz and comparator-norm adaptivity in online learning},
  author={Mhammedi, Z. and Koolen, W. M},
  booktitle={Conference on Learning Theory},
  pages={2858--2887},
  year={2020},
  organization={PMLR},
  url={http://proceedings.mlr.press/v125/mhammedi20a/mhammedi20a.pdf}
}

@InProceedings{FlaspohlerOCMOOM21,
  title = 	 {Online Learning with Optimism and Delay},
  author =       {Flaspohler, G. E. and Orabona, F. and Cohen, J. and Mouatadid, S. and Oprescu, M. and Orenstein, P. and Mackey, L.},
  booktitle = 	 {Proc. of the 38th International Conference on Machine Learning},
  pages = 	 {3363--3373},
  year = 	 {2021},
  editor = 	 {Meila, M. and Zhang, T.},
  volume = 	 {139},
  series = 	 {Proceedings of Machine Learning Research},
  month = 	 {18--24 Jul},
  publisher =    {PMLR},
  url = 	 {https://proceedings.mlr.press/v139/flaspohler21a.html},
}

@inproceedings{TammelinBJB15,
  title={Solving heads-up limit {Texas} hold'em},
  author={Tammelin, O. and Burch, N. and Johanson, M. and Bowling, M.},
  booktitle={Twenty-fourth international joint conference on artificial intelligence},
  year={2015},
  url = {https://poker.cs.ualberta.ca/publications/2015-ijcai-cfrplus.pdf}
}

@article{Blackwell56,
  title={An analog of the minimax theorem for vector payoffs},
  author={Blackwell, D.},
  journal={Pacific Journal of Mathematics},
  volume={6},
  number={1},
  pages={1--8},
  year={1956},
}

@article{HartMC00,
  title={A simple adaptive procedure leading to correlated equilibrium},
  author={Hart, S. and Mas-Colell, A.},
  journal={Econometrica},
  volume={68},
  number={5},
  pages={1127--1150},
  year={2000},
  publisher={Wiley Online Library},
  url={https://www.ma.imperial.ac.uk/~dturaev/Hart0.pdf}
}

@article{OrabonaP21,
  title={Parameter-free Stochastic Optimization of Variationally Coherent Functions},
  author = {Orabona, F. and P\'al, D.},
  journal={arXiv preprint arXiv:2102.00236},
  year=2021,
  url={https://arxiv.org/abs/2102.00236}
}

@article{grand2024solving,
  title={Solving optimization problems with {Blackwell} approachability},
  author={Grand-Cl{\'e}ment, Julien and Kroer, Christian},
  journal={Mathematics of Operations Research},
  volume={49},
  number={2},
  pages={697--728},
  year={2024},
  publisher={INFORMS}
}

@techreport{gordon2005no,
  title={No-regret algorithms for structured prediction problems},
  author={Gordon, Geoffrey J},
  year={2005},
  institution={Carnegie Mellon University. Center for Automated Learning and Discovery}
}

@book{brinkhuis2020convex,
  title={Convex analysis for optimization},
  author={Brinkhuis, Jan},
  year={2020},
  publisher={Springer}
}

@inproceedings{blackwell1954controlled,
  title={Controlled random walks},
  author={Blackwell, David},
  booktitle={Proceedings of the international congress of mathematicians},
  volume={3},
  pages={336--338},
  year={1954}
}

@phdthesis{Gordon99,
  title={Approximate solutions to Markov decision processes},
  author={Gordon, G. J.},
  year={1999},
  school={Carnegie Mellon University},
  url={https://www.proquest.com/docview/304499958}
}

@inproceedings{ZhangCP22,
  title={{PDE}-based optimal strategy for unconstrained online learning},
  author={Zhang, Z. and Cutkosky, A. and Paschalidis, I.},
  booktitle={International Conference on Machine Learning},
  pages={26085--26115},
  year={2022},
  organization={PMLR},
  url={https://proceedings.mlr.press/v162/zhang22d.html}
}

@inproceedings{CarmonH22,
  title={Making {SGD} parameter-free},
  author={Carmon, Y. and Hinder, O.},
  booktitle={Conference on Learning Theory},
  pages={2360--2389},
  year={2022},
  organization={PMLR},
  url={https://proceedings.mlr.press/v178/carmon22a.html}
}

@inproceedings{IvgiHC23,
  title={{DoG} is {SGD}'s Best Friend: A Parameter-Free Dynamic Step Size Schedule},
  author={Ivgi, M. and Hinder, O. and Carmon, Y.},
  booktitle={International Conference on Machine Learning},
  year={2023},
  url={https://arxiv.org/abs/2302.12022}
}

@inproceedings{ChzhenGS23,
  title={Parameter-free projected gradient descent},
  author={Chzhen, E. and Giraud, C. and Stoltz, G.},
  booktitle={International Conference on Machine Learning},
  url={https://arxiv.org/abs/2305.19605},
  year={2023}
}

@inproceedings{FarinaKS21,
  title={Faster game solving via predictive {Blackwell} approachability: Connecting regret matching and mirror descent},
  author={Farina, G. and Kroer, C. and Sandholm, T.},
  booktitle={Proceedings of the AAAI Conference on Artificial Intelligence},
  volume={35},
  number={6},
  pages={5363--5371},
  year={2021},
  url={https://arxiv.org/pdf/2007.14358.pdf}
}

@inproceedings{ChenKLO22,
  title={Better parameter-free stochastic optimization with {ODE} updates for coin-betting},
  author={Chen, K. and Langford, J. and Orabona, F.},
  booktitle={Proceedings of the AAAI Conference on Artificial Intelligence},
  volume={36},
  number={6},
  pages={6239--6247},
  year={2022},
  url={https://ojs.aaai.org/index.php/AAAI/article/view/20573}
}

@inproceedings{GrandClementK21,
  title={Conic {Blackwell} {A}lgorithm: Parameter-Free Convex-Concave Saddle-Point Solving},
  author={Grand-Cl\'{e}ment, J. and Kroer, C.},
  booktitle={Advances in Neural Information Processing Systems},
  volume={34},
  pages={9587--9599},
  publisher={Curran Associates, Inc.},
  year={2021},
  url={https://proceedings.neurips.cc/paper_files/paper/2021/file/4f87658ef0de194413056248a00ce009-Paper.pdf}
}

@inproceedings{ChakrabartiGK24,
  title={Extensive-Form Game Solving via {Blackwell} Approachability on Treeplexes},
  author={Chakrabarti, D. and Grand-Cl\'{e}ment, J. and Kroer, C.},
  booktitle={Advances in Neural Information Processing Systems},
  volume={37},
  pages={35257--35287},
  publisher={Curran Associates, Inc.},
  year={2024},
  doi={10.52202/079017-1111},
  url={https://proceedings.neurips.cc/paper_files/paper/2024/file/3e2aeb66481dd63a32421bf032b70384-Paper-Conference.pdf}
}

@inproceedings{CutkoskyM24,
  title={Fully Unconstrained Online Learning},
  author={Cutkosky, A. and Mhammedi, Z.},
  booktitle={Advances in Neural Information Processing Systems},
  volume={37},
  pages={10148--10201},
  publisher={Curran Associates, Inc.},
  year={2024},
  doi={10.52202/079017-0326},
  url={https://proceedings.neurips.cc/paper_files/paper/2024/file/13cd22c32c1330decd69c13cf8cadc0a-Paper-Conference.pdf}
}

@article{Tammelin2014,
    Author        = {Tammelin, Oskari},
    Title         = {Solving Large Imperfect Information Games Using {CFR+}},
    Journal       = {arXiv preprint arXiv:1407.5042},
    url           = {https://arxiv.org/abs/1407.5042},
    Year          = {2014},
}

@unpublished{Kwon26approachability,
  author = {Joon Kwon},
  title = {A general class of adaptive algorithms for approachability, regret and optimization},
  year = {2026},
  note = {In preparation},
}

@phdthesis{kwon2016mirror,
  title={Mirror descent strategies for regret minimization and approachability},
  author={Kwon, Joon},
  year=2016,
  school={Universit{\'e} Pierre-et-Marie-Curie}
}

@article{kwon2021refined,
  title={Refined approachability algorithms and application to regret minimization with global costs.},
  author={Kwon, Joon},
  journal={Journal of Machine Learning Research},
  volume={22},
  pages={200--1},
  year={2021}
}

@article{perchet2014approachability,
title={Approachability, regret and calibration: {I}mplications and equivalences},
journal={Journal of Dynamics and Games},
volume ={1},
number ={2},
pages ={181--254},
year ={2014},
author ={Perchet, V.}
}

@inproceedings{CutkoskyDM23,
  title={Mechanic: A learning rate tuner},
  author={Cutkosky, A. and Defazio, A. and Mehta, H.},
  booktitle={Advances in neural information processing systems},
  volume={36},
  pages={47828--47848},
  year={2023}
}

\appendix
\section{Proofs of the elementary results in Section \ref{sec:notation}}\label{app:prf-preliminaries}
\begin{proof}[Proof of Lemma \ref{lem:squared-norm-subgradient}]
Fix $\bd\in\partial\phi(\by)$. Apply the subgradient inequality to $(1+t)\by$ and let $t\to0$ from both sides to obtain $\langle\bd,\by\rangle=\|\by\|^{2}$.
For every $\bv$ and $t>0$, the subgradient and triangle inequalities give
\[
 \langle\bd,\bv\rangle
 \leqslant\frac{\phi(\by+t\bv)-\phi(\by)}{t}
 \leqslant\|\by\|\|\bv\|+\frac{t}{2}\|\bv\|^{2}\,.
\]
Taking the limit as $t\to0^{+}$ proves $\|\bd\|_{\star}\leqslant\|\by\|$, including at $\by=\bzero$.
The inner-product identity gives the reverse inequality when $\by\ne\bzero$.
\end{proof}

\begin{proof}[Proof of Lemma \ref{lem:product-convexity}]
Adding the two strong-convexity inequalities proves strong convexity with respect to the product norm, since the weighted squared distances sum to its square.
By definition of the dual norms,
\[
\sup_{\|(\bx,\by)\|\leqslant1}\left\{\langle\bg,\bx\rangle+\langle\bz,\by\rangle\right\}
= \sup_{\substack{r,s\geqslant0\\r^{2}+s^{2}\leqslant1}}\left\{\frac{r\|\bg\|_{\mathcal{X},\star}}{\sqrt{\mu}}+\frac{s\|\bz\|_{\mathcal{Y},\star}}{\sqrt{\nu}}\right\}
=\sqrt{\frac{\|\bg\|_{\mathcal{X},\star}^{2}}{\mu}+\frac{\|\bz\|_{\mathcal{Y},\star}^{2}}{\nu}}\,.
\]
The second equality follows from Cauchy--Schwarz, with equality when $(r,s)$ is the unit vector proportional to $(\|\bg\|_{\mathcal{X},\star}/\sqrt{\mu},\|\bz\|_{\mathcal{Y},\star}/\sqrt{\nu})$.
If both dual vectors are zero, both sides are zero.
\end{proof}

\begin{proof}[Proof of Lemma \ref{lem:minimum-growth}]
Optimality gives $\bzero\in\partial(f+\iota_{\mathcal{D}})(\bx^{\star})$.
Continuity of $r$ allows the subdifferential sum rule, which supplies $\bxi\in\partial r(\bx^{\star})$ with $-\bxi\in\partial(p+\iota_{\mathcal{D}})(\bx^{\star})$.
For every $\by\in\mathcal{D}$, the strong-convexity inequality gives
\[
 p(\by)-p(\bx^{\star})
 \geqslant-\langle\bxi,\by-\bx^{\star}\rangle+\frac{\mu}{2}|\by-\bx^{\star}|^{2}\,.
\]
Adding $r(\by)-r(\bx^{\star})$ proves \eqref{eq:minimum-growth}.
The subgradient characterization of strong convexity in \citet[Lemma~4.2]{Orabona19} follows by taking a directional limit along the segment from $\bx^{\star}$ to $\by$; the same proof applies to a seminorm in a real normed space.
\end{proof}
\begin{proof}[Proof of Lemma \ref{lem:logarithm-bound}]
For $x<1$, the identity and bound
\[
-\ln(1-x)-x
=x^{2}\int_{0}^{1} \! \frac{r}{1-rx}\,\mathrm{d}r
\leqslant\frac{x^{2}}{2(1-\max\{x,0\})}
\]
prove the first two inequalities.
For $x>0$, the function $\ln x-1+1/x$ is minimized at $x=1$, proving \eqref{eq:logarithm-ratio-bounds}.
\end{proof}

\section{The \texorpdfstring{\ac{FTRL}}{FTRL} regret equality}
\label{app:ftrl-equality}

The following decomposition is due to \citet[Lemma~7.1]{Orabona19}. The algebraic proof below applies to arbitrary real vector spaces, including Hilbert spaces.

\begin{theorem}[\ac{FTRL} regret equality]
\label{thm:ftrl-equality}
Let $\mathcal{X}$ be a real vector space, let $\mathcal{X}^{\star}$ denote its algebraic dual, and let $\mathcal{V}\subseteq\mathcal{X}$ be nonempty.
For $t=1,\dots,T+1$, let $\psi_{t}\colon\mathcal{X}\to\mathbb{R}\cup\{+\infty\}$ be regularizers, and for $t=1,\dots,T$, let $\bg_{t}\in\mathcal{X}^{\star}$.
Define
\[
 F_{t}(\bx)
 \coloneqq\psi_{t}(\bx)+\sum_{i=1}^{t-1}\langle \bg_{i},\bx\rangle
 \quad\text{and}\quad
 \bx_{t}\in\argmin_{\bx\in\mathcal{V}}\ F_{t}(\bx)\,.
\]
Assume these minima are attained, and all terms below are finite.
Then, for every comparator $\bu\in\mathcal{X}$ for which these terms are defined,
\begin{equation}
\label{eq:ftrl-regret-equality}
\begin{aligned}
 \sum_{t=1}^{T} \langle \bg_{t}, \bx_{t}-\bu \rangle
 &=\psi_{T+1}(\bu)-\psi_{1}(\bx_{1})+\sum_{t=1}^{T}\left(F_{t}(\bx_{t})-F_{t+1}(\bx_{t+1})+\langle \bg_{t},\bx_{t}\rangle\right)\\
 &\quad+F_{T+1}(\bx_{T+1})-F_{T+1}(\bu)\,.
\end{aligned}
\end{equation}
If $\bu\in\mathcal{V}$, the last difference is nonpositive.
\end{theorem}
\begin{proof}
The sum of consecutive objective values telescopes to
$F_{1}\left(\bx_{1}\right)-F_{T+1}\left(\bx_{T+1}\right)$.
Use $F_{1}=\psi_{1}$ and
$F_{T+1}\left(\bu\right)=\psi_{T+1}\left(\bu\right)+\sum_{t=1}^{T} \langle \bg_{t}, \bu \rangle$ in the right-hand side of \eqref{eq:ftrl-regret-equality}.
Optimality of $\bx_{T+1}$ gives the final assertion.
\end{proof}

We call the three contributions the \emph{initial condition}, the \emph{cost of unknown losses}, and the \emph{final condition}. The initial condition compares the comparator's terminal regularizer with the minimum initial regularizer. Each summand in the cost of
unknown losses is a \emph{per-round cost}: it measures the price of choosing $\bx_{t}$ before observing $\bg_{t}$, including the change in regularization. The final condition compares the terminal minimizer with the comparator and is nonpositive for a feasible comparator.
The identity needs neither convexity nor differentiability.

The identity extends to time-dependent constraints by adding $\iota_{\mathcal{D}_{t}}$ to the regularizer. In our application, $\mathcal{D}_{t+1}\subseteq\mathcal{D}_{t}$, so the next iterate is feasible for the current objective, as required by the local stability proof. A comparator in $\mathcal{D}_{T+1}$ has a nonpositive final condition. The terminal regularizer serves only to analyze rounds $1,\dots,T$ and does not affect their predictions.

\section{The cone $\cCV$ need not be closed}
\label{sec:open_cone}

For $\mathcal{V}=\mathbb{R}_{+}$, the lifted cone is $\mathcal{C}_{\mathcal{V}}=\{(y,s):y\geqslant0,\ s>0\}\cup\{(0,0)\}$.
It is not closed: for every integer $n\geqslant1$, $(1,1/n)=(1/n)(n,1)$ belongs to $\mathcal{C}_{\mathcal{V}}$, whereas its limit $(1,0)$ does not.

This example shows that Lemma~12(1) of \citet{abernethy2011blackwell} needs the additional assumption that the convex cone is closed: in general, the bipolar of a convex cone is its closure. Their subsequent lifted construction assumes a compact convex decision set.

\end{document}